\documentclass[a4paper, 11pt]{article}

\usepackage{amssymb,latexsym,amsmath}     
\usepackage{graphicx,color}
\usepackage{algorithm,algpseudocode,bm} 

\newtheorem{theorem}{Theorem}

\newtheorem{remark}{Remark}

\newcommand{\R}{\mathbb{R}}

\usepackage{subcaption}

\begin{document}

\title{Numerical Hopf--Lax formulae for\\ Hamilton--Jacobi equations on unstructured geometries}
\author{S. Cacace\footnote{Dipartimento di Matematica, ``Sapienza'' Universit\`a di Roma, P.le Aldo Moro, 1, Roma, Italy; {\tt simone.cacace@uniroma1.it}}, R. Ferretti\footnote{Dipartimento di Matematica e Fisica, Universit\`a Roma Tre, Largo S. Leonardo Murialdo, 1, Roma, Italy; {\tt roberto.ferretti@uniroma3.it}}, G. Tatafiore\footnote{Dipartimento di Matematica, ``Sapienza'' Universit\`a di Roma, P.le Aldo Moro, 1, Roma, Italy; {\tt giulia.tatafiore@uniroma1.it}}}
\maketitle

\begin{abstract}
We consider a scheme of Semi-Lagrangian (SL) type for the numerical solution of Hamilton--Jacobi (HJ) equations on unstructured triangular grids. As it is well known, SL schemes are not well suited for unstructured grids, due to the cost of the point location phase; this drawback is augmented by the need for repeated minimization. In this work, we consider an existing, monotone version of the scheme, that works only on the basis of node values, and adapt the algorithm to the case of an unstructured grid, using the connectivity information. Then, applying a quadratic refinement to the numerical solution, we improve accuracy at the price of some extra computational cost. The scheme can be applied to both time-dependent and stationary HJ equations; in the latter case, we also study the construction of a fast policy iteration solver and of a parallel version. We perform a theoretical analysis of the two versions, and validate them with an extensive set of examples, both in the time-dependent and in the stationary case.
\end{abstract}

\section{Introduction}

Born in the 1950s in the framework of environmental fluid dynamics, semi-Lagrangian (SL) schemes are typically used to approximate hyperbolic equations, or advection terms in more complex evolution operators, using characteristics in numerical form. The simplest case is that of the constant-coefficient advection equation
\[
\begin{cases}
u_t + a\cdot Du = 0 & (x,t)\in\R^d\times\R^+, \\
u(x,0) = u_0(x),
\end{cases}
\]
where the vector $a\in\R^d$ is the advection speed and $Du$ denotes the gradient of the function $u$.
In this case, the solution may be represented by the well-known formula
\begin{equation}\label{eq:char}
    u(x,t) = u_0(x-at),
\end{equation}
which is in turn discretized in SL form as
\begin{equation}\label{eq:schema_trasp}
v_j^{n+1} = I[V^n](x_j-a\Delta t).
\end{equation}
In \eqref{eq:schema_trasp}, $I[V](x)$ represents an interpolation operator based on the values $v_i$ of the vector $V$ (associated to a space grid of nodes $x_i$), computed at the point $x$. Moreover, given a time discretization of step $\Delta t$, the notation $v_j^n$ stands for an approximation of the solution at the node $x_j$ and time $t_n=n\Delta t$, and $V^n$ stands for the vector collecting all such values.

Our interest here is towards the Hamilton--Jacobi (HJ) equation
\begin{equation}\label{eq:hj}
\begin{cases}
u_t + H(Du) = 0 & (x,t)\in\R^d\times\R^+, \\
u(x,0) = u_0(x),
\end{cases}
\end{equation}
with a strictly convex Hamiltonian $H$. In this case, the characteristic-based representation formula \eqref{eq:char} admits a nonlinear generalization, namely the following {\em Hopf--Lax formula}
\begin{equation}\label{eq:hl}
    u(x,t)=\inf_{y\in\R^d}\left\{u_0(y)+tH^*\left(\frac{x-y}{t}\right)\right\},
\end{equation}
where $H^*$ is the Legendre transform of $H$, defined as
\[
  H^*(q)=\sup_{p\in\R^d}\{p\cdot q - H(p)\}.
\]
We recall that, if $H$ is {\em strongly coercive}, that is, if
$$
\lim_{|p|\to\infty} \frac {H(p)}{|p|} = +\infty,
$$
then $H^*$ is also well defined, strictly convex and strongly coercive.

The similarity between \eqref{eq:hl} and \eqref{eq:char} is more apparent if we set $a=(x-y)/t$, and \eqref{eq:hl} is recast in the form
\begin{equation}\label{eq:hl2}
    u(x,t)=\inf_{a\in\R^d}\{u_0(x-at)+tH^*(a)\},
\end{equation}
which has been widely used for the SL discretization of \eqref{eq:hj}, as well as of the Bellman equation of Dynamic Programming (see \cite{FF01,FF13} and the references therein).

The numerical discretization of \eqref{eq:hj} has been extensively studied; a general review on the topic is provided in \cite{FF16}. However, a relatively small part of this literature is devoted to the construction of efficient schemes on complex geometries, which require unstructured space grids. Among the relevant works in this direction, we quote in particular the researches on monotone, ENO--WENO and Discontinuous Galerkin schemes (see the review \cite{S}) and those using representation formulae of Hopf--Lax or Dynamic Programming type \cite{BR,KOQ,SV}.

This latter technique has the advantage of allowing the use of large time steps, and has therefore been used to construct schemes of semi-Lagrangian type.
For \eqref{eq:hl2}, the typical SL discretization has the structure
\begin{equation}\label{eq:sl}
    v_j^{n+1}=\min_{a\in\R^d}\{I[V^n](x-a\Delta t)+\Delta tH^*(a)\},
\end{equation}
and implies two steps, crucial for efficiency: interpolation and minimization. In this formulation of the scheme, efficiency may have a serious drop when using unstructured grids, due to the complexity of locating the point $x-a\Delta t$ on the space triangulation, a phase to be repeated at each computation of the function to be minimized.

In this work, we will discretize again the representation formula \eqref{eq:hl} to obtain an approximation scheme for \eqref{eq:hj}, but drop the interpolation phase and restrict the computation to nodes. We remark that this use of the Hopf--Lax formula has various precedents (for example, \cite{BFZ,GS}; however, our interest here is to improve efficiency on unstructured grids to which this work is ultimately directed. To this end, the basic idea will be to perform minimization with a walk on the space grid, much in the spirit of \cite{CF}.We will then adapt this technique to problems with Dirichlet boundary conditions and to the stationary case
\begin{equation}\label{eq:hj_staz}
\lambda u + H(Du) = f(x),
\end{equation}
focusing in particular on a fast solver based on Policy Iteration (PI).

The outline of the paper is as follows. In Section \ref{sec:constr_ev} we will construct the scheme in the evolutive case and prove its convergence on the basis of the Barles--Souganidis theory \cite{BS}, while Section \ref{sec:constr_staz} retraces the construction of the scheme for the stationary case \eqref{eq:hj_staz}. Sections \ref{sec:unstructured}--\ref{sec:refinement} treat respectively the unstructured implementation and a technique of quadratic refinement to improve accuracy of the scheme. Last, we present a numerical validation in Section \ref{sec:examples}, and draw some conclusions in Section \ref{sec:conclu}.

\section{Construction of the scheme, time-dependent case}\label{sec:constr_ev}

We start by listing the minimal assumptions on the problem \eqref{eq:hj}.

\subsubsection*{Basic assumptions}

\begin{itemize}

\item
$H(\cdot)$ is $C^2$ with bounded second derivatives, so that
\[
\max_p|D^2H(p)|_2 \le M_H;
\]
moreover, it is strictly convex and  strongly coercive;

\item
$u_0$ is bounded and Lipschitz continuous.

\end{itemize}

We recall that, under these assumptions, a unique solution of \eqref{eq:hj} exists in the viscosity sense. For any $t>0$, the solution is also semiconcave, that is, it satisfies a unilateral upper bound on the second incremental ratios (see \cite{FF13}).

In the numerical approximation, we start by assuming that the space grid covers the whole of $\R^d$. We will denote by $v_j^n \quad (j,n\in\mathbb N)$ the approximation of the solution at the node $x_j$ and time $t_n=n\Delta t$ (for $j\in\mathbb N$ and $n\ge 0$), and $V^n$ will stand for the infinite vector collecting all such values. We avoid any assumption on both the regularity of the grid and the shape of the elements (although in the numerical examples they will be taken as triangles). We will also denote by $\Delta x$ the maximum diameter of the grid elements.

\subsection{Restricting the Hopf--Lax formula to nodes}\label{subsec:nodes}

We start by rewriting the Hopf--Lax formula \eqref{eq:hl} at a node $x_j$ and on a single time step, from $t_n$ to $t_{n+1}$:
\begin{equation}\label{eq:hl_deltat}
    u(x_j,t_{n+1})=\inf_{y\in\R^d}\left\{u(y,t_n)+\Delta tH^*\left(\frac{x_j-y}{\Delta t}\right)\right\}.
\end{equation}
We note that, due to the strong coercivity of $H$ and boundedness of $u_0$, the $\inf$ operator can be replaced by a $\min$. Indeed, since the function
$$
\mathcal F_j^0(y,t)= u_0(y)+ t H^* \left(\frac {x_j-y}{t}\right)
$$
is continuous and coercive, there exists a minimum point $\bar y$ for every $x\in\mathbb R^d,\, t\in \mathbb R^+$. 
Moreover, given the basic assumptions above, the search for the minimum point can be restricted to a compact set, as proved in \cite{FF01}.

Then, we restrict the computation to nodes, obtaining the fully discrete version
\begin{equation}\label{eq:hl_discr}
\begin{cases}
    \displaystyle v_j^{n+1}=\min_{k\in\mathbb{N}}\left\{v_k^n+\Delta tH^*\left(\frac{x_j-x_k}{\Delta t}\right)\right\}, \\
    v_j^0=u_0(x_j).
\end{cases}
\end{equation}
In previous works (see, e.g., \cite{FF01}), the discretization of \eqref{eq:hl_deltat} has been carried out by taking $y\in\R^d$ and using a descent-type minimization algorithm, with an interpolated value for $u(y,t_n)$. In this case, at the price of some drop in accuracy, we avoid both phases of interpolation and iterative minimization, which would result in a high computational cost. 

For the moment, we will neglect boundary conditions, and analyse convergence of the scheme assuming an infinite grid covering the whole of $\R^d$, whereas in a second phase we will examine the case of a bounded domain. Here, we will work in the framework of Barles--Souganidis theory \cite{BS}. As it is well-known, this convergence theory requires the scheme to be monotone, $L^\infty$-stable and consistent. In our case, the scheme is invariant for the addition of constants, since, for any $c\in\R$,
\[
  \min_{k\in\mathbb{N}}\left\{(v_k^n+c)+\Delta tH^*\left(\frac{x_j-x_k}{\Delta t}\right)\right\} = \min_{k\in\mathbb{N}}\left\{v_k^n+\Delta tH^*\left(\frac{x_j-x_k}{\Delta t}\right)\right\}+c,
\]
so that $L^\infty$ stability may be deduced from monotonicity via the results in \cite{CT}.

\subsubsection*{Monotonicity and $L^\infty$ stability}

Monotonicity of the scheme has already been proved elsewhere, but we include the proof for completeness. Assume that $V$ and $W$ are two vectors of nodal values such that $V\ge W$ element by element. We want to prove that
\begin{equation}\label{eq:monot}
    \min_{k\in\mathbb{N}}\left\{v_k+\Delta tH^*\left(\frac{x_j-x_k}{\Delta t}\right)\right\} \ge \min_{m\in\mathbb{N}}\left\{w_m+\Delta tH^*\left(\frac{x_j-x_m}{\Delta t}\right)\right\}.
\end{equation}
Let the indices achieving the minimum on the left- and right-hand side of \eqref{eq:monot} be denoted, respectively, by $k^*$ and $m^*$. Then, \eqref{eq:monot} may be rewritten as
\begin{equation*}
    v_{k^*}+\Delta tH^*\left(\frac{x_j-x_{k^*}}{\Delta t}\right) \ge w_{m^*}+\Delta tH^*\left(\frac{x_j-x_{m^*}}{\Delta t}\right).
\end{equation*}
On the other hand, $m^*$ achieves the minimum, and this implies that
\begin{equation}\label{eq:monot2}
    w_{k^*}+\Delta tH^*\left(\frac{x_j-x_{k^*}}{\Delta t}\right) \ge w_{m^*}+\Delta tH^*\left(\frac{x_j-x_{m^*}}{\Delta t}\right),
\end{equation}
but, since $V\ge W$, we also have that
\begin{equation*}
    v_{k^*}+\Delta tH^*\left(\frac{x_j-x_{k^*}}{\Delta t}\right) \ge w_{k^*}+\Delta tH^*\left(\frac{x_j-x_{k^*}}{\Delta t}\right).
\end{equation*}
Then, using \eqref{eq:monot2}, we obtain a fortiori that \eqref{eq:monot} is satisfied. The scheme \eqref{eq:hl_discr} is therefore monotone and $L^\infty$ stable.

\subsubsection*{Consistency}

Let now $u$ be a smooth solution of \eqref{eq:hj}. In its usual sense, consistency of the scheme requires that for all $j$, as $\Delta x,\Delta t\to 0$,
\begin{equation}\label{eq:cons}
    L(x_j,\Delta x,\Delta t) = \frac{1}{\Delta t}\left(u(x_j,t_{n+1})-\min_{k\in\mathbb{N}}\left\{u(x_k,t_n)+\Delta tH^*\left(\frac{x_j-x_k}{\Delta t}\right)\right\}\right) \to 0.
\end{equation}
Rewriting now $u(x_j,t_{n+1})$ via \eqref{eq:hl_deltat}, we get
\begin{eqnarray}\label{eq:cons2}
    L(x_j,\Delta x,\Delta t) & = & \frac{1}{\Delta t}\left(\min_{y\in\R^d}\left\{u(y,t_n)+\Delta tH^*\left(\frac{x_j-y}{\Delta t}\right)\right\}\right. \nonumber \\
    && \left. -\min_{k\in\mathbb{N}}\left\{u(x_k,t_n)+\Delta tH^*\left(\frac{x_j-x_k}{\Delta t}\right)\right\}\right).
\end{eqnarray}
We therefore have to estimate the difference between the minimum of the smooth function
\begin{equation}\label{eq:F_j}
\mathcal F_j^n(y,\Delta t) = u(y,t_n)+\Delta tH^*\left(\frac{x_j-y}{\Delta t}\right)
\end{equation}
and its minimum restricted to a grid with step $\Delta x$; for simplicity, we will carry out the computation in one space dimension and use the shorthand notation $\mathcal F_j(y)$ for \eqref{eq:F_j}.

Since at the minimum point $\bar y$ we have $\mathcal F_j'(\bar y)=0$, the difference $\mathcal F_j(x_k)-\mathcal F_j(\bar y)$, for $x_k-\bar y=O(\Delta x)$, may be estimated via a Taylor expansion centred at $\bar y$ as
\[
\mathcal F_j(x_k)-\mathcal F_j(\bar y) = \frac{\mathcal F_j''(\eta)}{2}(x_k-\bar y)^2 \le \frac{\|\mathcal F_j''\|_\infty}{2}\Delta x^2
\]

This requires to give an upper bound on the second derivative of $\mathcal F_j$. We have
\[
\mathcal F_j''(y) = u_{xx}(y,t_n)+ \Delta t\frac{1}{\Delta t^2}{H^*}''\left(\frac{x_j-y}{\Delta t}\right).
\]
The first term in the right-hand side is bounded for a smooth solution, while the second, given the boundedness of ${H^*}''$, is of order $1/\Delta t$. Then, the difference between the minimum of $\mathcal F_j$ and its minimum value at the nodes is $O(\|{H^*}''\|_\infty\Delta x^2)$, and the consistency error reads
\[
L(x_j,\Delta x,\Delta t) = \frac{1}{\Delta t}O\left(\frac{\Delta x^2}{\Delta t}\right) = O\left(\frac{\Delta x^2}{\Delta t^2}\right).
\]
Therefore, the scheme is consistent under an inverse CFL condition $\Delta x=o(\Delta t)$ (this latter consistency condition also appears in \cite{BFZ,GS}). For example, taking the relationship $\Delta t\sim\Delta x^{1/2}$, the scheme turns out to be first-order w.r.t. $\Delta x$.

\begin{remark}
The previous study of the local truncation error aims at characterizing the accuracy of the scheme, while \cite{BS} adopts a weaker concept of consistency. More precisely, the Barles--Souganidis theorem would require that, for any function $\phi\in C^\infty(\R^d\times[0,T])$,
\begin{equation}\label{eq:BS_consistency}
\lim_{\Delta t\to 0}\frac{\displaystyle\phi(x_j,t_{n+1})-\min_k\left\{\phi(x_k,t_n)+\Delta tH^*\left(\frac{x_j-x_k}{\Delta t}\right)\right\}}{\Delta t} = \phi_t(x_j,t_n)+H(D\phi(x_j,t_n)).
\end{equation}
However, we show that \eqref{eq:BS_consistency} also follows from the previous arguments. In fact, writing $H$ as the Legendre transform of $H^*$, we have
\begin{eqnarray*}
\phi_t(x_j,t_n)+H(D\phi(x_j,t_n)) & = & \phi_t(x_j,t_n)+\max_{q\in\R^d}\left\{q\cdot D\phi(x_j,t_n)-H^*(q)\right\} \\
& = & \frac{\phi(x_j,t_{n+1})-\phi(x_j,t_n)}{\Delta t} \\
&& +\max_y\left\{\frac{\phi(x_j,t_n)-\phi(y,t_n)}{\Delta t}-H^*\left(\frac{x_j-y}{\Delta t}\right)\right\}+O(\Delta t),
\end{eqnarray*}
where we have and approximated $\phi_t$ via the time incremental ratio, the directional derivative $q\cdot D\phi$ with the incremental ratio with step $\Delta t$ in the direction
$$
q = \frac{x_j-y}{\Delta t},
$$
and collected all the error terms in a $O(\Delta t)$. Simplifying now the terms $\pm\phi(x_j,t_n)/\Delta t$ and rewriting the $\max$ as a $\min$, we finally get
\begin{equation}\label{eq:BS_consistency_2}
\phi_t(x_j,t_n)+H(D\phi(x_j,t_n)) = \frac{\displaystyle\phi(x_j,t_{n+1})-\min_y\left\{\phi(y,t_n)+\Delta tH^*\left(\frac{x_j-y}{\Delta t}\right)\right\}}{\Delta t}+O(\Delta t).
\end{equation}
Using \eqref{eq:BS_consistency_2} in the right-hand side of \eqref{eq:BS_consistency}, consistency in the sense of \cite{BS} finally comes down to comparing $\mathcal F_j^n(y,\Delta t)$ with $\mathcal F_j^n(x_k,\Delta t)$, and may be obtained by the arguments already used above.
\end{remark}

We can finally summarize this convergence analysis in the following

\begin{theorem}\label{thm:CONV}
Let the basic assumptions hold, and assume, in addition, that $\Delta x=o(\Delta t)$. Then, the numerical solution $V^n$ defined by \eqref{eq:hl_discr} converges to the solution $u(x,t_n)$ of \eqref{eq:hj} locally uniformly on $\R^d\times [0,T]$ as $\Delta t,\Delta x\to 0$ for any $T>0$.
\end{theorem}

\subsection{Error estimates}

In order to obtain explicit error estimates for the numerical solution $V^n$, we recall that the exact solution $u$ of \eqref{eq:hj} is Lipschitz continuous (with, say, a Lipschitz constant $L_u$), and that the speed of propagation $(x-y)/\Delta t$ achieving the minimum in \eqref{eq:hl} is bounded, i.e., that we can restrict the argument of $H^*$ to a bounded set, in which $H^*$ is also Lipschitz continuous, with a Lipschitz constant $L_{H^*}$. Then, comparing \eqref{eq:hl_deltat} and \eqref{eq:hl_discr} and denoting respectively by $y_j$ and $k_j$ the two minimizers, we have a first one-sided bound of the form
\begin{eqnarray}\label{eq:stima1}
    u(x_j,t_{n+1}) - v_j^{n+1} & = & u(y_j,t_n)+\Delta tH^*\left(\frac{x_j-y_j}{\Delta t}\right) - v_{k_j}^n-\Delta tH^*\left(\frac{x_j-x_{k_j}}{\Delta t}\right) \nonumber \\
    & \le & u(x_{k_j},t_n)+\Delta tH^*\left(\frac{x_j-x_{k_j}}{\Delta t}\right) - v_{k_j}^n-\Delta tH^*\left(\frac{x_j-x_{k_j}}{\Delta t}\right) \nonumber \\
    & \le & u(x_{k_j},t_n) - v_{k_j}^n,
\end{eqnarray}
where we have used the fact that $x_{k_j}$ is not in general a minimizer in \eqref{eq:hl_deltat}. Iterating back \eqref{eq:stima1} and taking into account that $v_j^0=u_0(x_j)$, we obtain
\[  
u(x_j,t_{n+1}) - v_j^{n+1} \le 0.
\]
As for the reverse estimate, we denote by $x_{k_j^*}$ the grid node closest to $y_j$, so that $|x_{k_j^*}-y_j|\le\Delta x$. Then, we get:
\begin{eqnarray}\label{eq:stima2}
    v_j^{n+1} - u(x_j,t_{n+1}) & = & v_{k_j}^n+\Delta tH^*\left(\frac{x_j-x_{k_j}}{\Delta t}\right) - u(y_j,t_n)-\Delta tH^*\left(\frac{x_j-y_j}{\Delta t}\right)  \nonumber \\
    & \le & v_{k_j^*}^n+\Delta tH^*\left(\frac{x_j-x_{k_j^*}}{\Delta t}\right) - u(y_j,t_n)-\Delta tH^*\left(\frac{x_j-y_j}{\Delta t}\right)  \nonumber \\
    & = & \left(v_{k_j^*}^n - u(x_{k_j^*},t_n)\right) + \left(u(x_{k_j^*},t_n) - u(y_j,t_n)\right) \nonumber \\
    && + \Delta t\left|H^*\left(\frac{x_j-x_{k_j^*}}{\Delta t}\right) - H^*\left(\frac{x_j-y_j}{\Delta t}\right)\right| \nonumber \\
    & = & v_{k_j^*}^n - u(x_{k_j^*},t_n) + L_u\Delta x + \Delta t L_{H^*}\frac{\Delta x}{\Delta t}.
\end{eqnarray}
Finally, denoting by $U^n$ the vector of node values of the exact solution at $t_n$, iterating back \eqref{eq:stima2} for $t_n\le T$ and including \eqref{eq:stima1}, we obtain the estimate:
\[  
0 \le V^n - U^n \le T(L_u+L_{H^*})\frac{\Delta x}{\Delta t},
\]
where the inequalities are clearly understood component-wise.

\subsection{Boundary conditions}\label{subsec:boundary}

In case Dirichlet conditions are enforced in \eqref{eq:hj}, we assume that the HJ equation is posed on a bounded open set $\Omega\subset\R^d$, and that the boundary datum $b(x)$ is assigned at any boundary point $x\in\partial\Omega$, for a given continuous function $b$. Then, the HJ equation takes the form
\begin{equation}\label{eq:hj_boundary}
    \begin{cases}
        u_t + H(Du) = 0 & (x,t)\in\Omega\times\R^+, \\
        u(x) = b(x) & x\in\partial\Omega, \\
        u(x,0) = u_0(x), & x\in\Omega.
    \end{cases}
\end{equation}
We recall that the Dirichlet condition in \eqref{eq:hj_boundary} should be understood in the {\it weak sense} (see \cite{FF13} for details). In this case, up to a shift and change of sign of the time variable, the solution of \eqref{eq:hj_boundary} may be characterized, via dynamic programming arguments (see \cite{BCD}), as the value function of the optimal control problem
\begin{equation}\label{eq:cost}
    \min_{\alpha(\cdot)} \int_t^{T\wedge\tau_\alpha(x,t)} H^*(\alpha(s))ds + u_0(X_\alpha(T\wedge\tau_\alpha(x,t))),
\end{equation}
where $T>0$, $t\in[0,T]$, $\alpha\in L^1(\R^+,\R^d)$ and $X_\alpha(\cdot)$ satisfies
\begin{equation}\label{eq:dynamics}
    \begin{cases}
        \dot X_\alpha(s) = \alpha(s) & s\ge t, \\
        X_\alpha(t) = x\in\overline\Omega.
    \end{cases}
\end{equation}
Here, $\tau_\alpha(x,t)$ is the {\it exit time} associated to the initial condition $(x,t)$ and the control $\alpha$, defined as
\begin{equation}\label{eq:exit}
\tau_\alpha(x,t) := \inf\{s\ge t:X_\alpha(s)\in\partial\Omega\},
\end{equation}
and we use the notation $a\wedge b=\min\{a,b\}$. We also assume that the final cost $u_0$ is compatible with the boundary stopping cost $b$, i.e., that $u_0(x)=b(x)$ for $x\in\partial\Omega$.

Then, setting $\alpha(s)\equiv (x_j-y)/(\Delta t\wedge \tau)$, for $x\in\overline\Omega$ the Hopf--Lax formula in the form \eqref{eq:hl_deltat} is replaced by
\begin{equation}\label{eq:hl_boundary1}
u(x_j,t_{n+1})=\min_{y\in\overline\Omega\atop\tau\in (0,\Delta t]}\mathcal F_j^n(y,\tau)
\end{equation}
where
\begin{equation}\label{eq:hl_boundary2}
\mathcal F_j^n(y,\tau) = \begin{cases}
u(y,t_n)+\Delta t H^*\left(\frac{x_j-y}{\Delta t}\right) & \text{ if } x_j,y\in\Omega \\[2pt]
b(y)+\tau H^*\left(\frac{x_j-y}{\tau}\right) & \text{ if } x_j\in\Omega,y\in\partial\Omega \\[2pt]
\min\left\{b(x_j),u(y,t_n)+\Delta t H^*\left(\frac{x_j-y}{\Delta t}\right)\right\} & \text{ if } x_j\in\partial\Omega,y\in\Omega \\[2pt]
\min\left\{b(x_j),b(y)+\tau H^*\left(\frac{x_j-y}{\tau}\right)\right\} & \text{ if } x_j,y\in\partial\Omega,
\end{cases}
\end{equation}
and the fully discrete scheme is finally obtained as
\begin{equation}\label{eq:scheme_boundary1}
v_j^{n+1}=\min_{x_k\in\overline\Omega\atop\tau\in (0,\Delta t]} F_j^n(x_k,\tau),
\end{equation}
in which $F_j^n$ is the fully discrete version of $\mathcal F_j^n$, defined in turn as
\begin{equation}\label{eq:scheme_boundary2}
F_j^n(x_k,\tau) = \begin{cases}
v_k^n+\Delta t H^*\left(\frac{x_j-x_k}{\Delta t}\right) & \text{ if } x_j,x_k\in\Omega \\[2pt]
b(x_k)+\tau H^*\left(\frac{x_j-x_k}{\tau}\right) & \text{ if } x_j\in\Omega,x_k\in\partial\Omega \\[2pt]
\min\left\{b(x_j),v_k^n+\Delta t H^*\left(\frac{x_j-x_k}{\Delta t}\right)\right\} & \text{ if } x_j\in\partial\Omega,x_k\in\Omega \\[2pt]
\min\left\{b(x_j),b(x_k)+\tau H^*\left(\frac{x_j-x_k}{\tau}\right)\right\} & \text{ if } x_j,x_k\in\partial\Omega.
\end{cases}
\end{equation}
The discretization \eqref{eq:scheme_boundary1}--\eqref{eq:scheme_boundary2} differs from \eqref{eq:hl_discr} only when $y=x_k\in\partial\Omega$ -- in this case, a minimization with respect to the exit time is also required. Therefore, for each $x_k$ on the boundary, we must compute
\[
\min_{\tau\le\Delta t}\mathcal F_j^n(x_k,\tau) = \min_{\tau\le\Delta t}\left\{b(x_k)+\tau H^*\left(\frac{x_j-x_k}{\tau}\right)\right\}.
\]
The term to be minimized is smooth and strictly convex, therefore its minimum is attained at the minimum between its only stationary point $\tau^*$ and $\Delta t$ (we will denote this optimal value as $\tau_{jk}$). We have:
\begin{equation}\label{eq:t_uscita}
\frac d {d\tau} \left\{v_k^n+\tau H^*\left(\frac{x_j-x_k}{\tau}\right)\right\} = H^*\left(\frac{x_j-x_k}{\tau}\right) - \frac{x_j-x_k}{\tau}\cdot\nabla H^*\left(\frac{x_j-x_k}{\tau}\right),
\end{equation}
which allows to find the stationary point, either explicitly or numerically, as soon as the analytic form of $H^*$ is known, via a restriction of $H^*$ to the direction $x_j-x_k$. In particular, in the numerical tests we will use the hamiltonian
\[
H(p) = a_0 + \frac 1 2 p^2,
\]
whose Legendre transform is
\[
H^*(q) = \frac 1 2 q^2 - a_0.
\]
Using this form of $H^*$, the stationary point in \eqref{eq:t_uscita} is computed, after some algebra, by solving the stationarity condition
\[
- a_0 - \frac 1 2 \frac{(x_j-x_k)^2}{\tau^2} = 0,
\]
which has no solution if $a_0>0$ (in which case the optimal stopping time is $\Delta t$), and otherwise has the solution
\[
\tau^* = \frac{|x_j-x_k|}{\sqrt{-2a_0}}.
\]
The time step to be used with the node $x_k$ is therefore
\[
\tau_{jk} = \begin{cases} \Delta t & \text{ if } a_0>0, \\ \min\left\{\Delta t, \frac{|x_j-x_k|}{\sqrt{-2a_0}} \right\} & \text{ otherwise.} \end{cases}
\]

\section{Construction of the scheme, stationary case}\label{sec:constr_staz}

We turn now to the stationary equation \eqref{eq:hj_staz}. In the time-marching framework, its solution can be seen as the asymptotic solution, for $t\to\infty$, of the time-dependent equation
\begin{equation}\label{eq:time_marching}
\begin{cases}
u_t + \lambda u + H(Du) = f(x) & (x,t)\in\R^d\times\R^+, \\
u(x,0) = u_0(x).
\end{cases}
\end{equation}
On the other hand, using Dynamic Programming arguments (see \cite{FF13}), it is possible to prove that the solution of \eqref{eq:hj_staz} (or, in other terms, the regime solution of \eqref{eq:time_marching}) satisfies
\begin{equation}\label{eq:dppio}
	u(x)=\inf_{\alpha(\cdot)}\left\{e^{-\lambda \Delta t}u(X_\alpha(\Delta t)) + \int_0^{\Delta t} \Big(f(X_\alpha(s))+H^*(\alpha(s))\Big)e^{-\lambda s}ds\right\},
\end{equation}
where $\alpha\in L^1(\R^+,\R^d)$ and $X_\alpha(s)$ solves the ODE system \eqref{eq:dynamics}.

In view of the $x$-dependence of $f$, \eqref{eq:dppio} cannot be given by an exact Hopf--Lax type form; therefore, setting $y=X_\alpha(\Delta t)$, and $\alpha(s)\equiv (x-y)/\Delta t$, and assuming smoothness as usual, we approximate the integral in \eqref{eq:dppio} as
\begin{equation}\label{eq:approx_integral}
\int_0^{\Delta t} \Big(f(y(s;\alpha))+H^*(\alpha(s))\Big)e^{-\lambda s}ds = \Delta t \left(f(x)+H^*\left(\frac{x-y}{\Delta t}\right)\right) + O\left(\Delta t^2\right).
\end{equation}
With this approximation, \eqref{eq:dppio} is in turn approximated as
\begin{equation}\label{eq:w}
w(x)=\min_{y\in\R^d}\left\{e^{-\lambda \Delta t}w(y) + \Delta t H^*\left(\frac{x-y}{\Delta t}\right) \right\} + \Delta t f(x).
\end{equation}
In a second discretization step, according to what has been done for the time-dependent case, we restrict both the computation and the minimum search to the set of grid nodes, obtaining
\begin{equation}\label{eq:hl_staz_approx}
v_j=\min_{k\in\mathbb{N}}\left\{e^{-\lambda \Delta t}v_k+\Delta tH^*\left(\frac{x_j-x_k}{\Delta t}\right)\right\}+\Delta t f(x_j).
\end{equation}
Equation \eqref{eq:hl_staz_approx} is a fixed-point equation, and it is easy to show that it has a contractive right-hand side (see \cite{FF13}); this ensures that it has a unique solution $V$. Concerning the convergence of $V$ to the solution of \eqref{eq:hj_staz}, it is possible to retrace the arguments in \cite{FF13}, along with those of the time-dependent case above, to obtain consistency in the sense of \cite{BS}, and a consistency estimate in the form
\begin{equation}\label{eq:staz_consistency1}
L(x_j,\Delta x,\Delta t) = \frac{1}{\Delta t}\left(O\left(\frac{\Delta x^2}{\Delta t}\right)+O\left(\Delta t^2\right)\right) = O\left(\frac{\Delta x^2}{\Delta t^2}\right)+O(\Delta t).
\end{equation}
We have then the following

\begin{theorem}
Let the basic assumptions hold, and assume in addition that $\Delta x=o(\Delta t)$. Then, the numerical solution $V^n$ defined by \eqref{eq:hl_staz_approx} converges to the solution $u(x,t_n)$ of \eqref{eq:hj_staz} locally uniformly on $\R^d$ as $\Delta t,\Delta x\to 0$.
\end{theorem}

While approximating the integral as in \eqref{eq:approx_integral} suffices to obtain a consistent scheme for $\Delta x=o(\Delta t)$, the term $O(\Delta t)$ becomes a bottleneck for the consistency rate. An easy computation shows that the truncation error \eqref{eq:staz_consistency1} is maximized taking $\Delta t\sim\Delta x^{2/3}$, and accordingly the consistency rate w.r.t. $\Delta x$ drops to 2/3. A way to increase the rate is to approximate the integral in the left-hand side of \eqref{eq:approx_integral} by a trapezoidal formula,
\begin{eqnarray}\label{eq:approx_integral2}
\int_0^{\Delta t} \Big(f(y(s;\alpha))+H^*(\alpha(s))\Big)e^{-\lambda s}ds & = & \frac{\Delta t} 2 \Big(f(x)+e^{-\lambda \Delta t}f(y) \nonumber \\
&&+\big(1+e^{-\lambda \Delta t}\big)H^*\Big(\frac{x-y}{\Delta t}\Big)\Big) + O\left(\Delta t^3\right)
\end{eqnarray}
(where we are still assuming $\alpha(s)\equiv (x-y)/\Delta t$ to be constant, and $f$ to be smooth enough),
which results in the scheme
\begin{equation}\label{eq:hl_staz_approx2}
v_j=\min_{k\in\mathbb{N}}\left\{e^{-\lambda \Delta t}v_k+\frac{\Delta t}{2}\left(\big(1+e^{-\lambda \Delta t}\big)H^*\left(\frac{x_j-x_k}{\Delta t}\right)+e^{-\lambda \Delta t} f(x_k)\right)\right\}+\frac{\Delta t}{2} f(x_j).
\end{equation}
In this case, given the increase of the approximation order for the integral of $f$, the consistency error reads
\begin{equation}\label{eq:staz_consistency2}
L(x_j,\Delta x,\Delta t) = \frac{1}{\Delta t}\left(O\left(\frac{\Delta x^2}{\Delta t}\right)+O\left(\Delta t^3\right)\right) = O\left(\frac{\Delta x^2}{\Delta t^2}\right)+O\left(\Delta t^2\right),
\end{equation}
the consistency rate is maximized under the relationship $\Delta t\sim\Delta x^{1/2}$ and its optimal value is 1.

Finally, we remark that, mixing the arguments of the time-dependent case with those of \cite[Chapter 8]{FF13} it is possible to prove error estimates of the form
\[
\| V^n - U^n\|_\infty \le C_1\Delta t^\gamma + C_2\frac{\Delta x}{\Delta t},
\]
with $\gamma=1$ for the scheme \eqref{eq:hl_staz_approx}, and $\gamma=2$ for the scheme \eqref{eq:hl_staz_approx2}.

\subsection{Boundary conditions}

Combining now the construction above with the techniques of the time-dependent case, we can modify the scheme for the case of weak Dirichlet conditions by adding a stopping cost on the boundary and considering the exit time as defined in \eqref{eq:exit}.
The form \eqref{eq:hl_staz_approx} of the scheme is then replaced by
\begin{equation}\label{eq:hl_boundary_staz1}
v_j=\min_{x_k\in\overline\Omega\atop\tau\in (0,\Delta t]} F_j(V,x_k,\tau)
\end{equation}
where $V=(v_1,\ldots,v_N)$  and $F_j$ is defined as
\begin{equation}\label{eq:hl_boundary_staz2}
F_j(V,x_k,\tau) = \begin{cases}
e^{-\lambda\Delta t}v_k+\Delta t H^*\left(\frac{x_j-x_k}{\Delta t}\right)+\Delta tf(x_j) & \text{ if } x_j,x_k\in\Omega \\[2pt]
e^{-\lambda\tau}b(x_k)+\tau H^*\left(\frac{x_j-x_k}{\tau}\right)+\tau f(x_j) & \text{ if } x_j\in\Omega,x_k\in\partial\Omega \\[2pt]
\min\left\{b(x_j),e^{-\lambda\Delta t}v_k+\Delta t H^*\left(\frac{x_j-x_k}{\Delta t}\right)+\Delta tf(x_j)\right\} & \text{ if } x_j\in\partial\Omega,x_k\in\Omega \\[2pt]
\min\left\{b(x_j),e^{-\lambda\tau}b(x_k)+\tau H^*\left(\frac{x_j-y}{\tau}\right)+\tau f(x_j)\right\} & \text{ if } x_j,x_k\in\partial\Omega.
\end{cases}
\end{equation}
The corresponding definition for the discretization \eqref{eq:hl_staz_approx2} may be obtained with obvious changes. The optimal value $\tau_{jk}$ of $\tau$ for $x_k\in\partial\Omega$ may be computed following the guidelines of the time-dependent case, although in this case we do not expect an explicit solution, even for quadratic Hamiltonians. On the other hand, all the $\tau_{jk}$ can be computed once and for all at the start of the iterative solution of \eqref{eq:hl_boundary_staz1}.

\subsection{Value iteration}

In the easiest approach, \eqref{eq:hl_staz_approx} (or \eqref{eq:hl_boundary_staz1}) may be solved by the so-called {\em value iteration} (VI), that is, via the iteration
\begin{equation}\label{eq:value_it}
v_j^{(n+1)}=\min_{x_k\in\overline\Omega\atop\tau\in (0,\Delta t]} F_j\left(V^{(n)},x_k,\tau\right)
\end{equation}
where the index $n$ should now be understood as an iteration number (an obvious modification applies for treating \eqref{eq:hl_staz_approx2}).
Since the right-hand side of \eqref{eq:hl_staz_approx} is a contraction, the iteration \eqref{eq:value_it} converges to its unique fixed point. The contraction coefficient, however, is given by $e^{-\lambda \Delta t}$, and the convergence becomes very slow as the approximation is refined. This supports the use of faster solvers for the fixed-point equation, in particular the so-called {\em Policy Iteration} (PI), which will be briefly reviewed here. Note that, for shortness, we will outline the basic ideas on the case of an infinite grid, the introduction of boundary conditions being straightforward.

\subsection{Exact policy iteration}

For a problem with $N$ unknowns (in our case, with a grid of $N$ nodes), the standard matrix form of the Policy Iteration algorithm reads
\begin{equation}\label{eq:PI}
\min_{\alpha \in U^N}\big(B(\alpha)V - c(\alpha)\big)=0,
\end{equation}
in which $B$ is a square $N\times N$ matrix, $V\in\R^N$ is the vector of the discretized solution, and $c\in\R^N$ is the vector of costs. After setting $n=0$ and choosing an initial policy $\alpha^{(0)}$, the PI evolves by alternating the two phases

\begin{enumerate}

\item {\em Policy evaluation}: $V^{(n)}$ is computed as the solution of the linear system
\[
B\left(\alpha^{(n)}\right)V^{(n)} = c\left(\alpha^{(n)}\right),
\]

\item {\em Policy improvement}: a new policy $\alpha^{(n+1)}$ is set as
\[
\alpha^{(n+1)}\in\arg\min_{\alpha \in U^N}\left(B(\alpha)V^{(n)} - c(\alpha)\right),
\]

\end{enumerate}
until a suitable convergence condition is satisfied.

Among the wide literature on policy iteration, we quote here the pioneering theoretical analysis of Puterman and Brumelle \cite{PB}, which have shown that the linearization procedure underlying policy iteration is conceptually equivalent to a Newton-type iterative solver. More recently, this result has been further generalized in \cite{BMZ09} in the following superlinear convergence result: 

\begin{theorem}[\cite{BMZ09}]\label{thm:PI}
Assume that:

\begin{itemize}

\item[\bf{(P1)}] The functions $B: U^N\to \R^{N\times N}$ and $c:U^N\to \R^N$ are continuous.

\item[\bf{(P2)}] For every $\alpha\in U^N$, the matrix $B(\alpha)$ is monotone, that is, $B(\alpha)$ is invertible and $B^{-1}(\alpha)\geq 0$ elementwise.

\end{itemize}

Then, there exists a unique $V$ solution of \eqref{eq:PI}. Moreover, the sequence $V^{(n)}$ generated by the PI algorithm satisfies the following:
\begin{enumerate}
\item $V^{(n+1)}\leq V^{(n)}$ for every $n\geq 0$;
\item $V^{(n)}\to V$ as $n\to+\infty$, for every initial policy $\alpha_0$;
\item $\|V^{(n+1)}-V\|=o\left(\|V^{(n)}-V\|\right)$ as $n\to+\infty$.
\end{enumerate}
\end{theorem}

In the case of the discretization \eqref{eq:hl_staz_approx}, we can recast the scheme in the form
\begin{equation}\label{eq:PI_1}
    \min_{k_j\in\mathbb{N}}\left\{-v_j+e^{-\lambda \Delta t}v_{k_j}+\Delta t\left(H^*\left(\frac{x_j-x_{k_j}}{\Delta t}\right)+f(x_j)\right)\right\}=0\quad (j=1,\ldots,N),
\end{equation}
where we understand the collection of all the indices $k_j$, i.e., the vector
\[
{\bm k}=(k_1,\ldots,k_N)^T
\]
as the policy.

To interpret \eqref{eq:PI_1} in the form \eqref{eq:PI}, it suffices to define the elements of the matrix $B(\bm k)$ and the vector $c(\bm k)$ as
\begin{equation}
    b_{jm}(\bm k)=\begin{cases}
        1 & \text{ if } m=j\ne k_j \\
        -e^{-\lambda \Delta t} & \text{ if } m=k_j\ne j \\
        1-e^{-\lambda \Delta t} & \text{ if } m=j=k_j \\
        0 & \text{ otherwise};
    \end{cases}
    \hspace{1.5cm}
    c_j(\bm k)=\Delta t\left(H^*\left(\frac{x_j-x_{k_j}}{\Delta t}\right)+f(x_j)\right).
\end{equation}
In this case, denoting by $k_j^{(n)}$ the $j$-th node index in the $n$-th policy $\bm k^{(n)}$, and changing all the signs, the policy evaluation phase is rewritten as the system
\begin{equation}\label{eq:PI_2}
    v_j^{(n)}-e^{-\lambda \Delta t}v_{k_j^{(n)}}^{(n)} = \Delta t\left(H^*\left(\frac{x_j-x_{k_j^{(n)}}}{\Delta t}\right)+f(x_j)\right) \quad (j=1,\ldots,N),
\end{equation}
while the policy improvement is computed as
\[
k_j^{(n+1)} \in \arg\min_{k_j\in\mathbb{N}}\left\{-v_j^{(n)}+e^{-\lambda \Delta t}v_{k_j}^{(n)}+\Delta t\left(H^*\left(\frac{x_j-x_{k_j}}{\Delta t}\right)+f(x_j)\right)\right\}=0\quad (j=1,\ldots,N).
\]

With these definitions, assumptions {\bf (P1)}--{\bf (P2)} of Theorem \ref{thm:PI} are satisfied. In fact, the continuity of $B$ and $c$ is obvious once the set of node indices is naturally endowed with the discrete topology. As for {\bf (P2)}, we note that an equivalent definition of monotone matrix requires the off-diagonal elements to be nonpositive (which is clearly the case) and the eigenvalues to have nonnegative real part; this latter property is true, as a consequence of Gershgorin's theorem. Note that monotonicity of the matrices $B$ implies their nonsingularity, and thus the unique solvability of the policy evaluation phase.

Clearly, the form \eqref{eq:hl_staz_approx2} of the scheme can be treated by changing the definition of the vector $c(\bm k)$ into
\[
c_j(\bm k)=\frac{\Delta t}2\left(\big(1+e^{-\lambda \Delta t}\big)H^*\left(\frac{x_j-x_{k_j}}{\Delta t}\right)+f(x_j)+e^{-\lambda \Delta t}f(x_{k_j})\right),
\]
and the policy evaluation step is changed accordingly.

\subsection{Modified policy iteration}

An adaptation of policy iteration to large sparse problems has been proposed as ``modified policy iteration'' (MPI) in \cite{PS}, and has also become a classical tool. In this case, the policy evaluation phase is carried out in iterative form, which may be convenient in very large problems, whenever a direct solver would be unfeasible. In the specific case of \eqref{eq:PI_1}, once given an initial guess $\bm k^{(0)}$ for the policy, the two nested cycles of MPI can be written as

\begin{enumerate}

\item {\em Policy evaluation}: $V^{(n)}$ is computed by iterating on $m$ the update formula
\[
    v_{j,m+1}^{(n)}=e^{-\lambda \Delta t}v_{k_{j}^{(n)},m}^{(n)}+\Delta t\left(H^*\left(\frac{x_j-x_{k_{j}^{(n)}}}{\Delta t}\right)+f(x_j)\right)
\]
up to a given stopping condition;

\item {\em Policy improvement}: a new policy $k_{j}^{(n+1)}$ is set as
\[
    k_{j}^{(n+1)}\in\arg\min_{k\in\mathbb{N}}\left\{-v_{j}^{(n)}+e^{-\lambda \Delta t}v_{k}^{(n)}+\Delta t\left(H^*\left(\frac{x_j-x_k}{\Delta t}\right)+f(x_j)\right)\right\}.
\]

\item {\em Stopping condition}: if $\bm{k}^{(n+1)}\ne\bm{k}^{(n)}$ go to 1, else stop.

\end{enumerate}

Here, we have denoted by $v_{j,m}^{(n)}$ the $j$-th element of the $m$-th iterate in the policy evaluation, by $v_{j}^{(n)}$ the corresponding final value of the iteration (i.e., the $j$-th element of $V^{(n)}$), and by $k_{j}^{(n)}$ the $j$-th element of the optimal policy at the $n$-th iteration.
Note that, since the policy is defined, in this scheme, by the set of indices $k_j$ of the nodes achieving the minima for $F_j$, we can consider the policy to be convergent as soon as none of these indices is any longer updated.

Finally, note that all three solvers (Value Iteration, Exact Policy Iteration, Modified Policy Iteration) necessarily converge to the same solution.

\section{Unstructured implementation}\label{sec:unstructured}
We provide some details for the actual implementation of the fully discrete scheme under consideration, on two-dimensional triangular grids. The aim is to compute the global minimum in \eqref{eq:hl_discr} as efficiently as possible, in terms of both space and time complexity, by means of walks along the grid nodes that locally decrease the right hand side of our scheme. To this end, once a triangulation (typically, Delaunay) is given, we only assume to work with the minimal data structures describing the domain geometry, namely vertex and triangle lists, plus a connectivity list which contains, for each triangle, the indices of its first triangle neighbours along its three edges. More precisely, we denote by $\mathcal{X}=\{x_j\}_{j=1,\dots,|\mathcal{X}|}$ the list of vertex coordinates $x_j\in\mathbb{R}^2$ where $|\mathcal{X}|$ is the total number of vertices, by $\mathcal{T}=\{T_j=(j_1,j_2,j_3)\}_{j=1,\dots,|\mathcal{T}|}$ the list of index triplets defining the $j-$th triangle with vertices $x_{j_1},x_{j_2},x_{j_3}$ where $|\mathcal{T}|$ is the total number of triangles, and by $\mathcal{
TN}=\{TN_j=(j_1,j_2,j_3)\}_{j=1,\dots,|\mathcal{T}|}$ the list of index triplets corresponding to the three neighbours $T_{j_1}, T_{j_2}, T_{j_3}, $ of the $j-$th triangle (for boundary elements one or two indices are set to zero if the corresponding neighbours are missing). Note that, in order to walk on the grid along the vertices, using $\mathcal{TN}$ is much cheaper than storing the first neighbours of each vertex, since the number of neighbours per vertex typically ranges from four up to eleven. Here we just store three indices per triangle, and we can use this connectivity list to move around a vertex, and accomplish the same task of finding vertex neighbours in constant time, by means of a custom function {\tt vertex\_neighbours}.

Now, the crucial point of our minimization procedure is the initialization of each walk. 
It is well known that the solutions of \eqref{eq:hj} can develop gradient discontinuities (shocks) that, while moving in time, can change completely the domain of dependence of the solution itself. At the discrete level this implies that the grid node achieving the global minimum in \eqref{eq:hl_discr} can suddenly change from one time step to the next, hence that some walk can be trapped in a local minimum. As a workaround to this scenario, we propose a simple but still effective procedure to handle nonsmooth solutions. The idea is to run, for each grid node $x_j$, different walks starting from a batch of neighboring nodes within a circle of radius proportional to $\Delta t$, in order to climb possible discontinuities and drive at least one walk to the correct minimum. In two dimensions, a reasonable compromise in terms of computational efforts is to use four points, obtained by displacing $x_j$ along coordinate axis, two per dimension in opposite directions, and accepting as correct the best of the values obtained by the four walks. More precisely, given a constant $C>0$, we consider the four directions $e_1=(1,0), e_2=(-1,0), e_3=(0,1), e_4=(0,-1)$ and the four displacements of each grid node $x_j$ for $j=1,\dots,|\mathcal{X}|$: 
 \begin{equation}
 x_{j,d}=x_j+C\Delta t\,e_d\qquad \mbox{for } d=1,\dots,4\,.
 \label{eq:displacements length}
 \end{equation}
 Denoting by 
 $$
 k_{j,d}=\displaystyle\arg\min_{k=1,\dots,|\mathcal{X}|}|x_{j,d}-x_k|
 $$
 the index of the closest vertex projection of $x_{j,d}$ on the grid, we build the index displacement lists $\mathcal{K}^d$ by setting $\mathcal{K}^d_j=k_{j,d}$ for $d=1,\dots,4$ and $j=1,\dots,|\mathcal{X}|$. 
 
 This yields the descent algorithm for the computation of the proposed node-restricted Hopf--Lax formula, presented in Algorithm \ref{ALG1}. 
\begin{algorithm}[!t]
\small
\begin{algorithmic}[1]
\label{alg:HL}
\State Given $\mathcal X$, $\mathcal T$, $\mathcal{TN}$, $\Delta t$, $N$, $u_0$, $H^*$, $b$
\State  $v_j^0\leftarrow u_0(x_j)$ for $j=1,\dots,|\mathcal{X}|$.
\For{$n=0:N$}
\For{$j=1:|\mathcal{X}|$}
\For{$d=1:4$}
\State WALK $\leftarrow$ true
\State $k_{\text{old}}$ $\leftarrow$ $\mathcal{K}^d_j$
 \While{WALK}
 \State $\mathcal{N}\leftarrow \{k_{\text{old}}\} \bigcup{\tt vertex\_neighbors}(k_{\text{old}})$
 \State $k_{\text{new}}\,\leftarrow\,\displaystyle\arg\min_{k\in \mathcal{N}}\mathcal{F}(k)\,,\qquad \mathcal{F}(k):=\min_{\tau\in (0,\Delta t]} F_j^n(x_k,\tau)$ \quad with $F_j^n$ given by \eqref{eq:scheme_boundary2}
\State WALK $\leftarrow$ $(k_{\text{new}}\neq k_{\text{old}})$
\State $k_{\text{old}}\leftarrow k_{\text{new}}$
\EndWhile
\State $v_{j,d}$ $\leftarrow$ $\mathcal{F}({k_{\text{new}}})$
\EndFor
\State $v_j^{n+1}\leftarrow \displaystyle\min_{d=1,\dots,4} v_{j,d}$
\EndFor
\EndFor 
\end{algorithmic}
\caption{{\tt Node\_Restricted\_HL}} \label{ALG1}
\end{algorithm}

We remark that the treatment of boundary conditions is embedded in the evaluation of $F_j^n$. Moreover, whenever the node $x_k$ lies on the boundary of the domain, the one-dimensional optimization of $F_j^n$ with respect to $\tau$ can be readily performed employing a standard root finder for its derivative, such as bisection or Newton method in case higher order derivatives can be explicitly computed.

As a rule of thumb for the choice of the constant $C$ in \eqref{eq:displacements length}, note that, in order to handle an abrupt change in the minimum caused by the transit of a singularity, this constant should be of the order of (or slightly greater than) the maximum speed of propagation of the solution. For a Hamiltonian function of the form $H(p)$ used here, it is well known that the local speed of propagation at $x$ is $|\nabla H(Du(x))|_2$; on the other hand, by a zeroth-order Taylor expansion,
\begin{eqnarray*}
|\nabla H(p)|_2 & = & \left|\nabla H(0) + D^2H(\bar p)p\right|_2 \\
& \leq & |\nabla H(0)|_2 + \max\left|D^2H(\bar p)\right|_2 \max |p|_2.
\end{eqnarray*}
Since the gradient of the solution of \eqref{eq:hj} has a nonincreasing norm, $\max |p|_2$ may be bounded by $L_0$, the Lipschitz constant of $u_0$ in the Euclidean norm, while $\left|D^2H(\bar p)\right|_2\le M_H$ by the basic assumptions. Then, a reasonable choice of $C$ would be
\[
C \gtrsim |\nabla H(0)|_2 + M_H L_0.
\]
For example, for the Hamiltonian $H(p)=\frac 1 2 |p|^2$ which will be used in the numerical tests, we have $\nabla H(0)=0$ and $\left|D^2H(\bar p)\right|_2\equiv 1$. In the stationary case, it is possible to replace $L_0$ with the Lipschitz constant of the time-discrete solution; for example, for the discrete formulation \eqref{eq:w} the Lipschitz constant of $w$ is
\[
L_w = \frac{L_f}{\lambda},
\]
with $L_f$ denoting the Lipschitz constant of $f$.

\begin{remark}
Note that, assuming that the solution is convex, we can completely drop the initialization of walks based on different seeds. Actually, in this case, the right hand side in \eqref{eq:hl_boundary1}
is also convex and has a unique global minimum. Then, inspired by \cite{CF}, we can further optimize the process by starting the walk from the node that attained the minimum in the previous time step, and this leads to a significant speed-up in the overall computation.
\end{remark}

We finally observe that, since we are interested in walks that strictly decrease $F_j^n$, namely the right hand side of \eqref{eq:scheme_boundary1}, 
the function {\tt vertex\_neighbors} can be combined with an index-based comparison to restrict the evaluation of $F_j^n$
only to the vertices that have not already been visited. This can clearly mitigate the overall computational cost of the algorithm in case the computation of $F_j^n$ 
is particularly expensive. 

The implementation of the proposed schemes for the stationary equation \eqref{eq:hj_staz} on triangular grids parallels the one for the evolutive case presented in Section \ref{sec:constr_ev}. In particular, we have to take into account the new terms related to the discount factor $\lambda$ and the source $f$ appearing in the right-hand side of \eqref{eq:hl_staz_approx}. Moreover, for the Value Iteration scheme \eqref{eq:value_it} we can readily apply Algorithm \ref{ALG1}, just reinterpreting the number of time steps as the number of fixed-point iterations, and adding a stopping criterion based on some norm of the difference between two successive iterations. On the other hand, for both Exact and Modified Policy Iteration schemes, we have to embed our vertex neighbours minimum search in the policy improvement phase. 

\section{Quadratic refinement of the minimum search}\label{sec:refinement}

One major drawback of the scheme in the form \eqref{eq:hl_discr} is the relatively low accuracy. We propose therefore a more accurate version, still working without a numerical minimization phase, which will be sketched for simplicity in two space dimensions. In this version, the refinement of the minimum search is performed via a local quadratic approximation, which improves accuracy of the scheme at the price of loosing its monotonicity.

\begin{figure}[!t]
\begin{center}
\includegraphics[height=4.5cm]{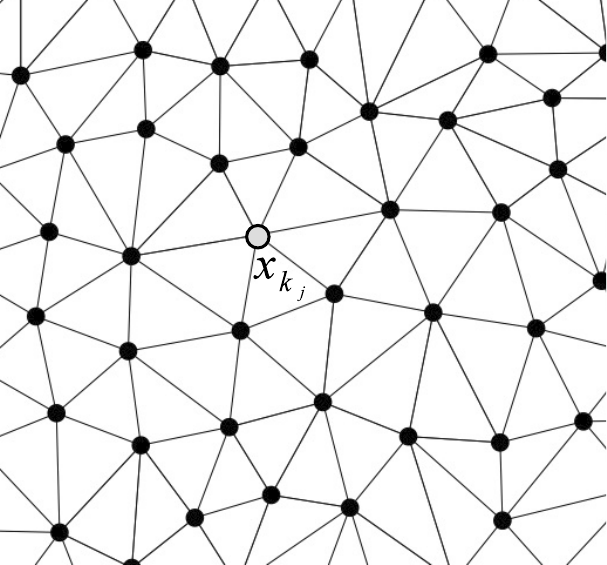} \hspace{1.5cm}
\includegraphics[height=4.5cm]{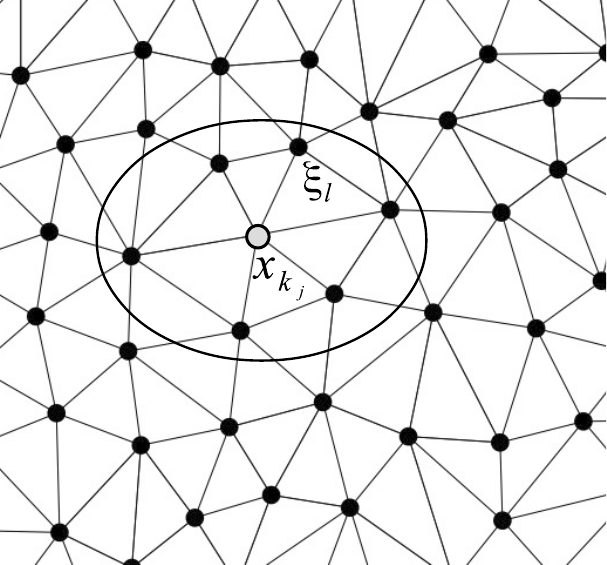}
\caption{Node achieving the minimum in \eqref{eq:hl_discr} (left) and related stencil for quadratic refinement (right).}\label{fig:stencil}
\end{center}
\end{figure}

Let $\mathcal F_j(y)$ be defined by \eqref{eq:F_j}, and denote by $k_j$ the index of the node achieving the minimum. Then, we consider the node $x_{k_j}$ along with its neighbours (as shown in Fig. \ref{fig:stencil}), and build a least squares quadratic approximation of $\mathcal F_j$ based on its values on such a stencil.

The construction of a quadratic approximation in the form
\begin{eqnarray}\label{eq:Q}
Q(y) & = & \frac 1 2 y^tAy - b^ty + c \nonumber \\
& = & \frac 1 2 a_{11} y_1^2 + a_{12}y_1y_2 + \frac 1 2 a_{22} y_2^2 - b_1y_1 - b_2y_2 + c,
\end{eqnarray}
with a symmetric Hessian matrix $A$, requires the six parameters $a_{11},a_{12},a_{22},b_1,b_2$ and $c$. Under the assumption that the triangulation is acute, the stencil contains at least six nodes, and this implies that the squared residual has a unique minimum (in particular, a zero residual if the number of nodes is precisely six, in which case $Q$ is an interpolant).
In \eqref{eq:Q} and in what follows, we neglect for simplicity of notation the dependence on the index $j$.

We briefly recall the adaptation of least squares approximation to the case under consideration. Let then $\xi_l=(\xi_{l,1},\xi_{l,2})$ for $l=1,\ldots,M$, denote a generic point (including $x_{k_j}$) in the stencil of $Q$, and $f_l$ denote the corresponding value of $\mathcal F_j(k)$, computed for $k$ such that $x_k=\xi_l$. Once defined the matrix $\Phi$ and the two vectors $a$ and $f$ as
\[
\Phi = \left(
\begin{matrix}
\displaystyle\frac 1 2 \xi_{1,1}^2 & \xi_{1,1}\xi_{1,2} & \displaystyle\frac 1 2 \xi_{1,2}^2 & -\xi_{1,1} & -\xi_{1,2} & 1 \\[3pt]
\vdots & \vdots & \vdots & \vdots & \vdots & \vdots \\[5pt]
\displaystyle\frac 1 2 \xi_{M,1}^2 & \xi_{M,1}\xi_{M,2} & \displaystyle\frac 1 2 \xi_{M,2}^2 & -\xi_{M,1} & -\xi_{M,2} & 1 \\
\end{matrix}
\right),\quad
a = \left(
\begin{matrix}
a_{11} \\
a_{12} \\
a_{22} \\
b_1 \\
b_2 \\
c \\
\end{matrix}
\right),\quad
f = \left(
\begin{matrix}
f_1 \\
\vdots \\
f_M \\
\end{matrix}
\right),
\]
the vector $a$ of unknown parameters of $Q(y)$ solves the system of normal equations
\begin{equation}\label{eq:eq_norm}
    \Phi^t\Phi a = \Phi^t f.
\end{equation}
Next, using the parameters obtained from \eqref{eq:eq_norm}, an easy computation gives for the minimum of $Q(y)$ the value
\[
    \min_{y\in\R^2} Q(y) = \frac{a_{22}b_1^2 - 2a_{12}b_1b_2 + a_{11}b_2^2}{2(a_{12}^2 - a_{11}a_{22})}+c.
\]
Finally, this value should replace $v_{k_j}^n+\Delta t H^*\left(\frac{x_j-x_{k_j}}{\Delta t}\right)$ in \eqref{eq:scheme_boundary2}, whenever $x_{k_j}\in\Omega$. If $x_{k_j}$ is on the boundary, then we expect that the stencil of quadratic refinement might not be unisolvent; in this case we take the node value instead of performing the refinement.

Note that the vector of unknown parameters $a$ may be written as $a=(\Phi^t\Phi)^{-1}\Phi^t f$, and the pseudoinverse $(\Phi^t\Phi)^{-1}\Phi^t$ depends on the local geometry of the grid, but not on the values in $f$, that is, it is a constant matrix associated to the node $x_{k_j}$. 

\begin{remark}\label{rem:acuteness}
Although we expect the stencil $\{\xi_l\}_{1\le l \le M}$ to be unisolvent (i.e., that $M\ge 6$) and the matrix $A$ to be positive definite, it might happen that one or both of these conditions fail to hold. For example, as we have just noted, unisolvency may be lost at the boundary; in our case, this may also occur at internal nodes since the Delaunay triangulation produced by the library {\tt Triangle} \cite{TRI} and used in the numerical tests is not necessarily acute. In these cases, the simplest countermeasure is to take the node value outright, without applying the quadratic refinement. Strategies to further improve the accuracy of quadratic refinement with Dirichlet conditions are outside the scope of this paper and will possibly be the object of a future work.
\end{remark}

\begin{remark}
While this technique proves useful in the time-dependent case (as we will show in the numerical examples), its application to the stationary equation is limited to the value iteration solver, which in some sense stems from a time-dependent formulation. On the contrary, it is not clear how this refinement could be introduced in a policy iteration setting. In order to compare the various solvers for the static problem, we will not apply quadratic refinement to this case.
\end{remark}




\subsubsection*{Consistency}

The consistency estimate for this version of the scheme will be derived by adapting the arguments used in \S \ref{subsec:nodes}, keeping the simplified 1-D setting and assuming in addition that $H$ has a bounded third derivative. First, note that, given a regular function $u$, the error obtained by a least squares, second-order polynomial approximation of $\mathcal F_j$ constructed on a fixed local stencils (of radius $O(\Delta x)$ and with a bounded number of nodes), retains the  same rate of the interpolation error, i.e., $O(\|\mathcal F_j'''\|_\infty\Delta x^3)$. We are not aware of results in this direction (this might be a somewhat classical result), but sketch here the arguments for completeness:

\begin{enumerate}
    \item Constructing a second order interpolating polynomial on any unisolvent subset of the stencil produces errors of order $O(\|\mathcal F_j'''\|_\infty\Delta x^3)$ on the other nodes of the stencil;

    \item None of these interpolators has a global residual lower that the least squares polynomial: the latter has therefore a squared residual of order (at worst) $O(\|\mathcal F_j'''\|_\infty\Delta x^3)$;

    \item Since the number of nodes of the stencil is bounded, the residual in the discrete $\infty$-norm is also $O(\|\mathcal F_j'''\|_\infty\Delta x^3)$;

    \item The least squares approximation appears then as an interpolation with perturbed data, where the perturbation is of the same order of the interpolation error. Then, the perturbed interpolation retains the same order of convergence of the unperturbed one.
    
\end{enumerate}

On the other hand,
\[
\mathcal F_j'''(y) = u_{xxx}(y,t_n) - \Delta t\frac{1}{\Delta t^3}{H^*}'''\left(\frac{x_j-y}{\Delta t}\right).
\]
Again, the first term in the right-hand side is bounded, while the second, assuming that ${H^*}'''$ is also bounded, is of order $1/\Delta t^2$. In this case, the difference between the minimum of $\mathcal F_j$ and the minimum of its interpolate is $O(\|{H^*}'''\|_\infty\Delta x^3)$, and the consistency error reads
\[
L(x_j,\Delta x,\Delta t) = \frac{1}{\Delta t}O\left(\frac{\Delta x^3}{\Delta t^2}\right) = O\left(\frac{\Delta x^3}{\Delta t^3}\right).
\]
Therefore, the scheme is again consistent under an inverse CFL condition $\Delta x=o(\Delta t)$, but with a higher rate for the same $\Delta t/\Delta x$ relationship. For example, the relationship $\Delta t\sim\Delta x^{1/2}$ which causes the monotone scheme to be consistent with order $O(\Delta x)$, results in this case in a consistency rate of $O(\Delta x^{3/2})$.

\section{Numerical examples}\label{sec:examples}


We present some numerical tests showing both the accuracy and the performance of the proposed schemes for evolutive and stationary HJ equations. Our implementation is written in C language, employing {\tt Triangle} \cite{TRI}, our favorite easy-to-use library for the generation of Delaunay triangular meshes in two dimensions, and {\tt SuiteSparse} \cite{SUS}, a sparse matrix library for the solution of linear systems associated to policy evaluation and quadratic refinement of the minimum search (although this choice might not be optimal in every respect).  
Unless otherwise stated, the numerical tests in serial mode are run on a MacBook Pro equipped with an Intel Core i5 dual-core with 2.9 GHz, 8 Gb RAM under OS 11.7.10, whereas the parallel codes are written in CUDA and running on a GPU server equipped with a Nvidia H100 GPU with 94Gb Ram and 14592 CUDA cores.



\subsection{Evolutive case}
\begin{figure}[!t]
    \centering
    \includegraphics[width=.49\textwidth]{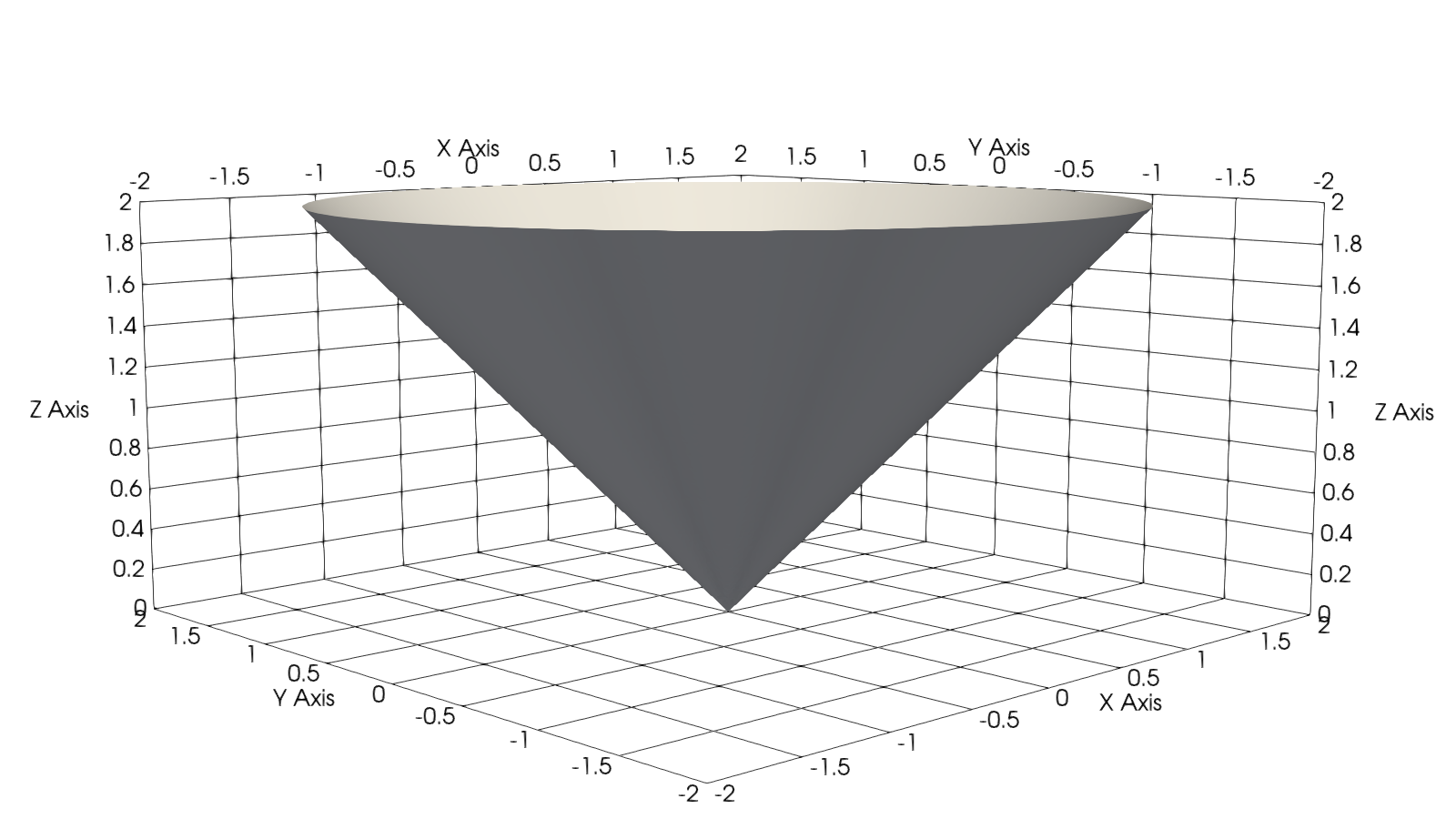}
    \includegraphics[width=.49\textwidth]{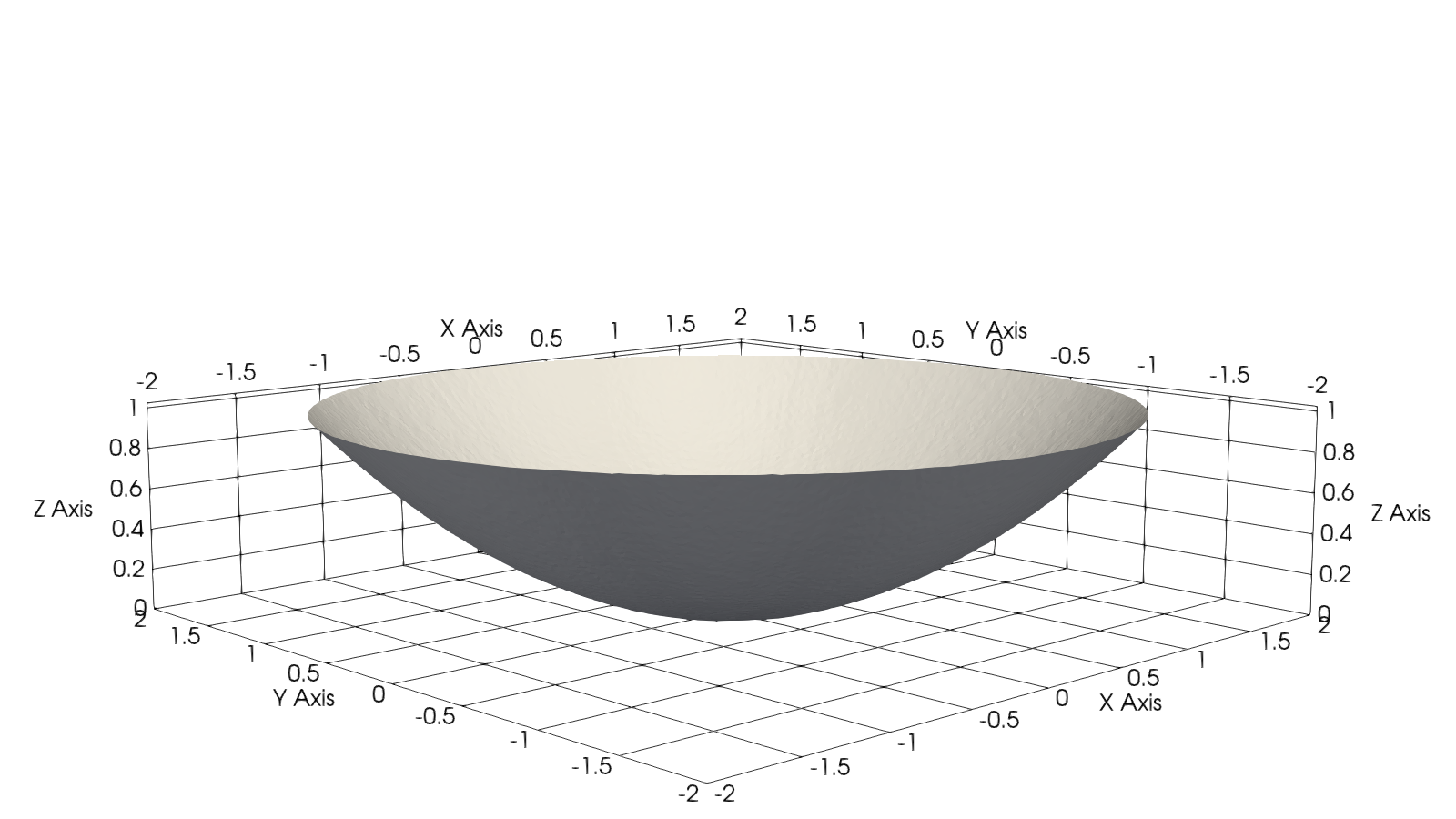}
    
    \caption{Test 1, numerical solution at initial time $t=0$ and final time $T=2$, for $\Delta x= 0.025$ and $\Delta t= 0.5 \sqrt{\Delta x}\approx 0.08$. }
    \label{fig:test1}
\end{figure}

\begin{figure}[!t]
    \centering
    \includegraphics[width=.49\textwidth]{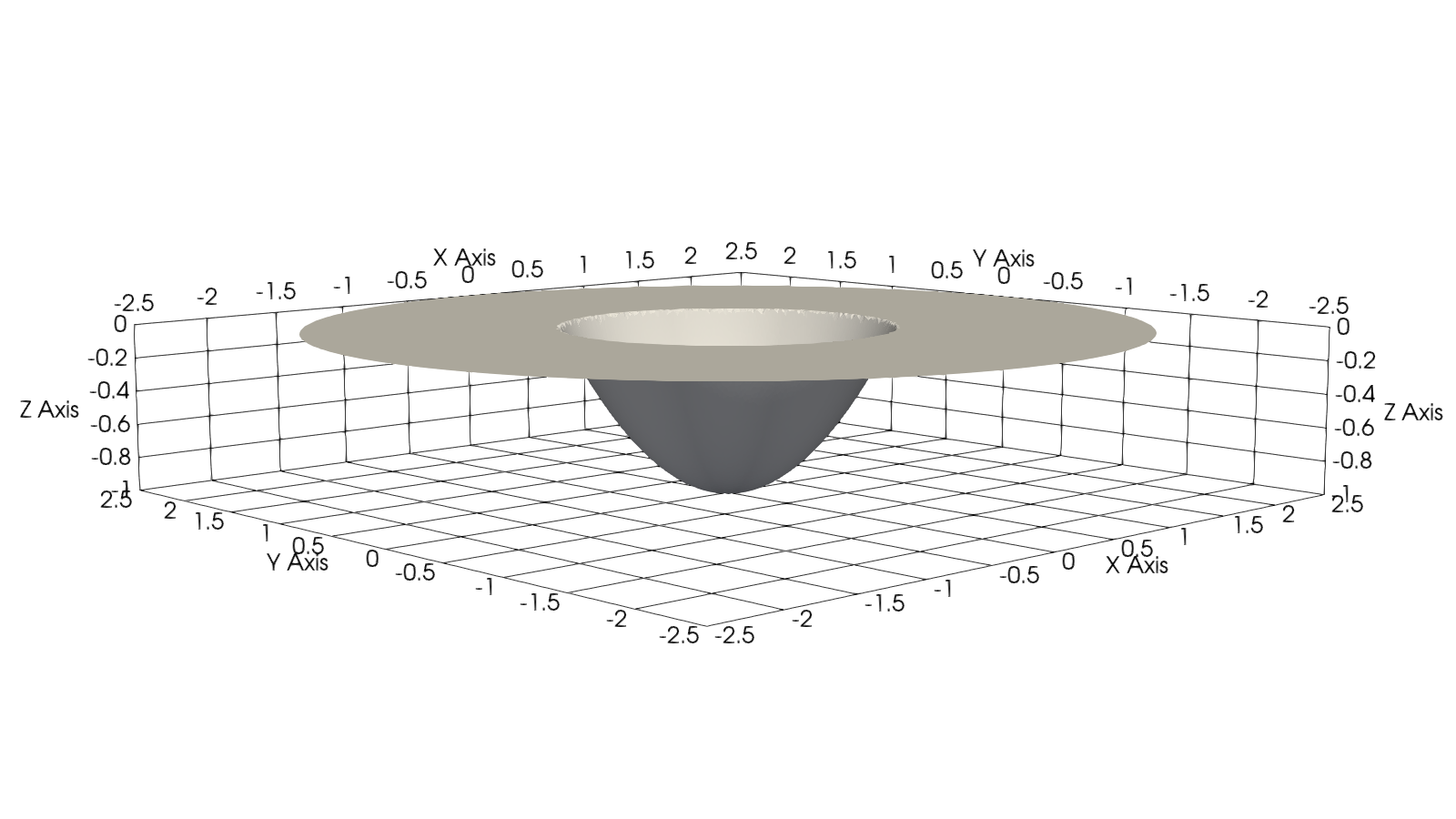}
    \includegraphics[width=.49\textwidth]{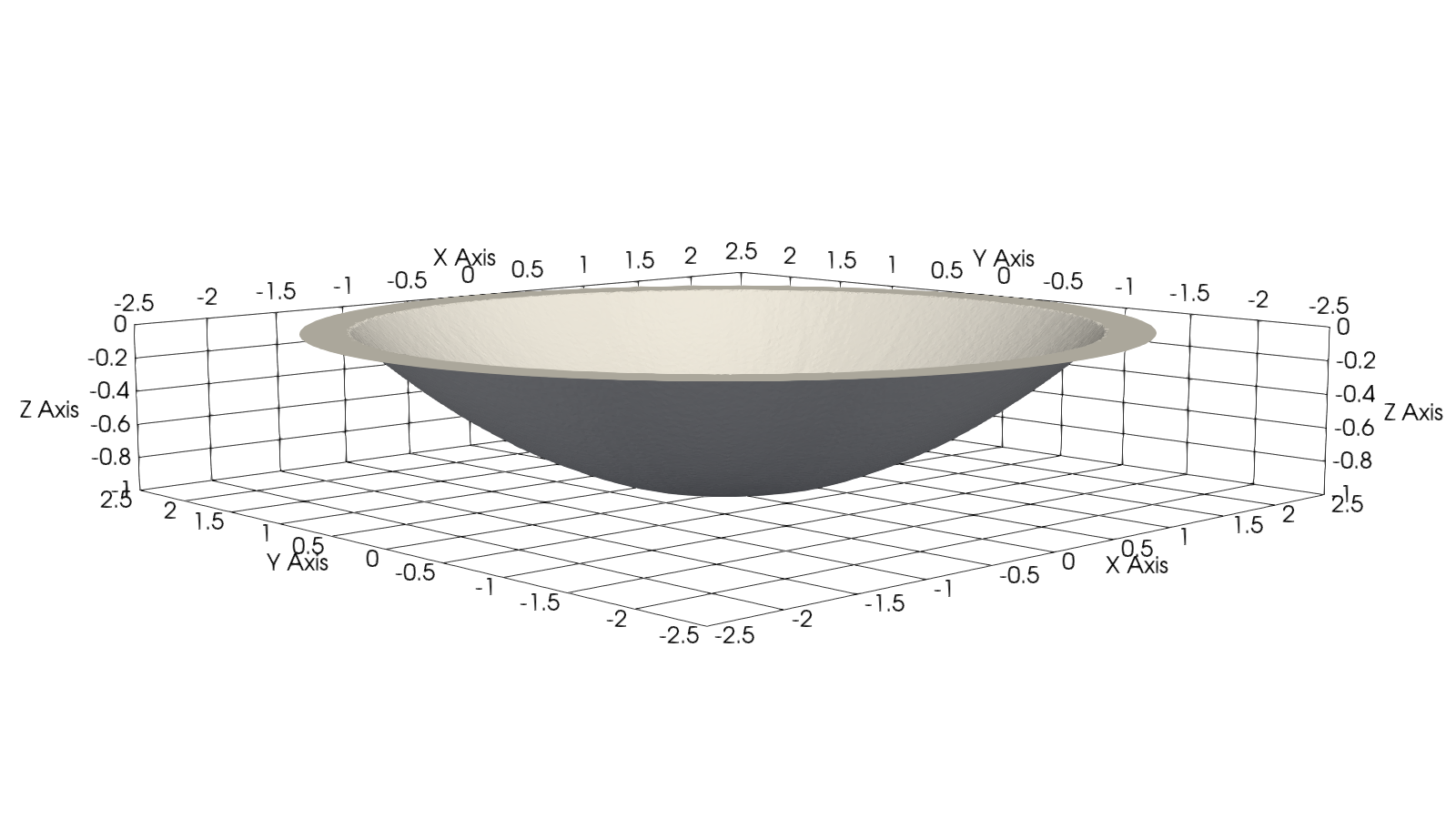}
    \caption{Test 2, numerical solution at initial time $t=0$ and final time $T=2$, for $\Delta x= 0.025$ and $\Delta t= 0.5 \sqrt{\Delta x}\approx 0.08$. }
    \label{fig:test2}
\end{figure}

For a first set of experiments in the time-dependent case, we consider the quadratic Hamiltonian $H:\mathbb{R}^2\to\mathbb{R}$ defined by $H(p)=\frac12|p|^2$, choose the two initial data 
$$
u_{0,1}(x)=|x|\,,\qquad u_{0,2}(x)=\min\left\{|x|^2-1\,,\,0\right\}\,,\qquad x\in\Omega\,,
$$
where the space domain $\Omega$ is a disc of $\mathbb{R}^2$, centered in the origin and with radius respectively equal to $2$ in the first case and $2.5$ in the second. Moreover, we fix the final time as $T=2$. In this setting, it is easily proved that $H$ and its Legendre transform coincide, and that the solutions of \eqref{eq:hj} with initial data $u_{0,1}$ and $u_{0,2}$ are given respectively by
$$
u_1(x,t)=\left\{
\begin{array}{ll}
     \displaystyle\frac{|x|^2}{2t} & \mbox{if } |x|\le t \\[10pt]
     |x|-\displaystyle\frac{t}{2} & \mbox{otherwise} 
\end{array}
\right.,\qquad
u_2(x,t)=\min\left\{ \frac{|x|^2}{2t+1} -1\,,\,0\right\}\,, 
$$
providing classical examples of convex and semiconcave solutions (see \cite{FF13}). In particular, the first solution exhibits an instantaneous regularization of the initial datum, whereas for the second one, the initial discontinuity of the gradient propagates along expanding circles. We observe that, for this final time $T$, the singularity remains in the interior of the domain. Initial conditions and final numerical solutions are shown in Figures \ref{fig:test1}--\ref{fig:test2}.

In order to perform an experimental analysis for the order of convergence of the scheme \eqref{eq:hl_discr}, we fulfill all the assumptions of Theorem \ref{thm:CONV}. In particular, choosing $\Delta t =c \sqrt{\Delta x}$, for a $c>0$, where $\Delta x$ is the maximal diameter of the elements of the triangulation, we ensure a consistency error of order $O(\Delta x)$. 
Then, denoting by $u_{\text{ex}}$ the exact solution (for $\text{ex}=1,2$), and by $t_N$ the final time $T$, we consider relative errors in the $L^\infty$ and $L^1$ norm,
\[
\begin{split}
&\mathcal{E}^\infty_{\Delta x}=\frac{\sup_{i}|u_{\text{ex}}(x_i,t_{N})-v_i^N|}{\sup_{i}|u_{\text{ex}}(x_i,t_{N})|},\\
&\mathcal{E}^1_{\Delta x}=\frac {\Delta x^2 \sum_i{|u_{\text{ex}}(x_i,t_{N})-v_i^N|}}{\Delta x^2 \sum_i{|u_{\text{ex}}(x_i,t_{N})|}}.
\end{split}
\]
Finally, we evaluate, under grid refinement, the experimental order of convergence given by $\text{EOC}\sim\log_2(\mathcal{E}_{\Delta x}/\mathcal{E}_{\Delta x/2})$, for errors in both $L^\infty$ and $L^1$ norms, denoted, respectively, by $\text{EOC}^\infty$ and $\text{EOC}^1$. 

The results for Test $1$ and Test $2$, associated with the initial conditions $u_{0,1}$ and $u_{0,2}$, are shown in Tables \ref{tab:Test1}--\ref{tab:Test2}, obtained with the coefficient $C$ of the displacements \eqref{eq:displacements length} set to $C=2$. In the tables we also report 
the sizes $|\mathcal X|$ and $|\mathcal T|$ of the vertex and triangle lists for each refinement of the grid.

Note that the experimental order of convergence follows the theoretical consistency analysis.
Moreover, test $2$ allows us to show the advantages of implementing the four displacement technique. In fact, in this case, there is a large part of the domain where the solution is constant and where the minimization walks would immediately stop at a local minimum if starting from the node $x_j$. The displacements, if their length $C\Delta t$ is large enough, allow the walks to bypass the singularity and find the right minimizing node.


\begin{table}[!t]
    \centering
    \begin{tabular}{|c|c|c|c|c|c|c|c|c|c|}
\hline
        $\Delta x$ & $\Delta t$ & $|\mathcal{X}|$ & $|\mathcal{T}|$  & $\mathcal{E}^1_{\Delta x}$ & $\text{EOC}^1$ & $\mathcal{E}^\infty_{\Delta x}$ & $\text{EOC}^\infty$\\
        \hline\hline
         0.1 & 0.1581 & 3142 & 6080 & 0.0523 & - &  0.0582 & -\\
\hline 
        0.05 & 0.1118 & 12360 & 24329 &  0.025 & 1.09 & 0.031 & 0.90 \\
\hline 
        0.025 & 0.0791 & 49077 & 97344 & 0.013 & 0.9404 &  0.0153 & 1.0225 \\
\hline  
        0.0125 & 0.0559 & 195420 & 389229 & 0.0060 &  1.0896 & 0.0068 & 1.1638 \\
\hline 
    \end{tabular}
    \caption{Test 1, experimental order of convergence for the basic version of scheme \eqref{eq:hl_discr}, with $\Delta t = 0.5 \sqrt{\Delta x}$.}
    \label{tab:Test1}

    \centering
    \begin{tabular}{|c|c|c|c|c|c|c|c|c|c|}
\hline
        $\Delta x$ & $\Delta t$ & $|\mathcal{X}|$ & $|\mathcal{T}|$ & $\mathcal{E}^1_{\Delta x}$ & $\text{EOC}^1$ & $\mathcal{E}^\infty_{\Delta x}$ & $\text{EOC}^\infty$\\
        \hline\hline
        0.1 & 0.1581 & 3142 & 6080 & 0.0918 & - & 0.0917 & - \\
\hline 
        0.05 & 0.1118 & 12360 & 24329 & 0.0415 & 1.1451 & 0.0435 & 1.0769\\
\hline 
        0.025 & 0.0791 & 49077 & 97344 & 0.0198 & 1.0675 &  0.0217 & 1.0024 \\
\hline  
        0.0125 & 0.0559 & 195420 & 389229 & 0.0094 &  1.0770 & 0.01050 & 1.0478 \\
\hline 
    \end{tabular}
    \caption{Test 2, experimental order of convergence for the basic version of scheme \eqref{eq:hl_discr}, with $\Delta t = 0.5 \sqrt{\Delta x}$.}
    \label{tab:Test2}
\end{table}

\subsection{Quadratic refinement}
The application of the quadratic refinement previously introduced is studied in the time-dependent case. We repeat test 1 and test 2 adding the minimum search on the quadratic approximation, and we collect the results in Tables \ref{tab:Test1_LS} and \ref{tab:Test2_LS}, respectively.

\begin{table}[!t]
    \centering
    \begin{tabular}{|c|c|c|c|c|c|c|c|c|c|}
\hline
        $\Delta x$ & $\Delta t$ & $|\mathcal{X}|$ & $|\mathcal{T}|$ &  $\mathcal{E}^1_{\Delta x}$ & $\text{EOC}^1$ & $\mathcal{E}^\infty_{\Delta x}$ & $\text{EOC}^\infty$\\
        \hline\hline
        0.1 & 0.1581 & 3142 & 6080 & 0.0112 & - & 0.0234 & - \\
\hline 
        0.05 & 0.1118 & 12360 & 24329 & 0.0038 & 1.5397 &  0.0088 &  1.4172 \\
\hline 
        0.025 & 0.0791 & 49077 & 97344 & 0.0025 &  0.6030 & 0.0059 & 0.5788 \\
\hline  
        0.0125 & 0.0559 & 195420 & 389229 & 0.0010 & 1.3899 & 0.0025 & 1.2247 \\
\hline 
    \end{tabular}
    \caption{Test 1, experimental order of convergence for scheme \eqref{eq:hl_discr} with quadratic refinement application, where $\Delta t = 0.5 \sqrt{\Delta x}$.}
    \label{tab:Test1_LS}

    \centering
    \begin{tabular}{|c|c|c|c|c|c|c|c|c|c|}
\hline
        $\Delta x$ & $\Delta t$ & $|\mathcal{X}|$ & $|\mathcal{T}|$ & $\mathcal{E}^1_{\Delta x}$ & $\text{EOC}^1$ & $\mathcal{E}^\infty_{\Delta x}$ & $\text{EOC}^\infty$ \\
        \hline\hline
        0.1 & 0.1581 & 3142 & 6080 & 0.0012 & - & 0.0118 & -\\
\hline 
        0.05 & 0.1118 & 12360 & 24329 & 0.0004 & 1.5727 &  0.0056  & 1.0812  \\
\hline 
        0.025 & 0.0791 & 49077 & 97344 & 0.0002 & 1.1114  & 0.0018 &  1.6059 \\
\hline  
        0.0125 & 0.0559 & 195420 & 389229  & 0.00008 & 1.1599 & 0.0014 & 0.3370  \\
\hline 
    \end{tabular}
    \caption{Test 2, experimental order of convergence for scheme \eqref{eq:hl_discr} with quadratic refinement application, where $\Delta t = 0.5 \sqrt{\Delta x}$.}    
    \label{tab:Test2_LS}
\end{table}

\begin{figure}[!t]
    \centering
    \includegraphics[width=0.49\textwidth]{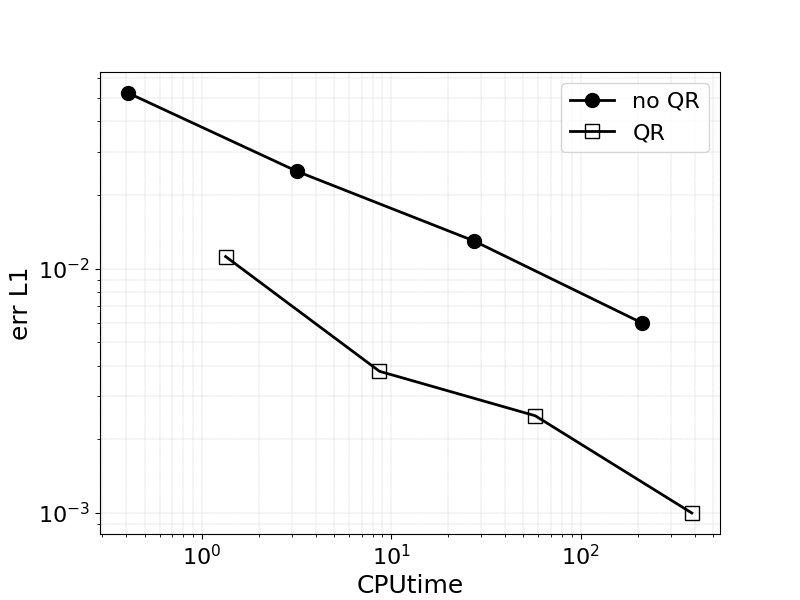}
    \includegraphics[width=0.49\textwidth]{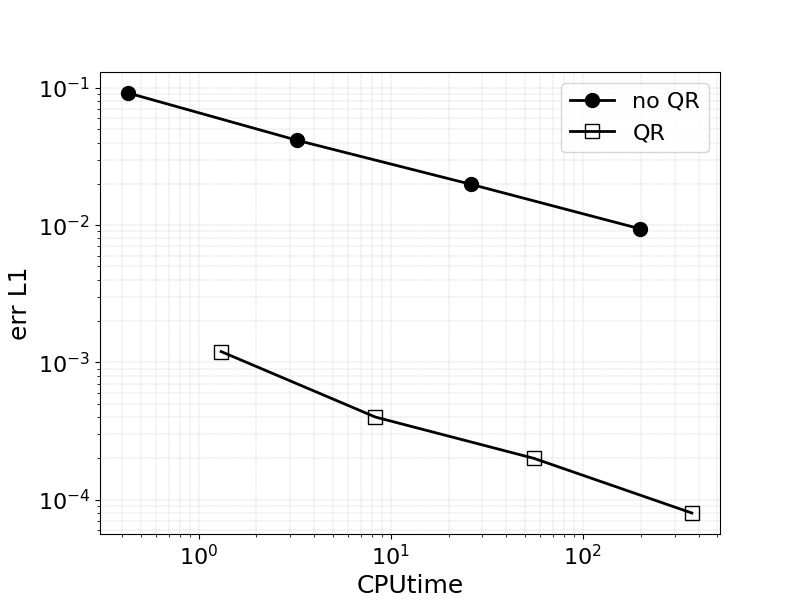}
    \caption{Test $1$ (left) and test $2$ (right) relative errors evolution, with respect to CPU times. Obtained applying scheme \eqref{eq:hl_discr} with and without quadratic refinement, respectively for $\Delta x= 0.1, 0.05, 0.025, 0.0125$.}
    \label{fig:QRcomparison}
\end{figure}

We observe that the implementation of this technique results in a relevant reduction of the relative errors, at the cost of an increased CPU time. Figure \ref{fig:QRcomparison} shows that the increase in accuracy is worth the higher complexity, by comparing the numerical errors in $L^1$-norm with respect to the CPU times for the scheme \eqref{eq:hl_discr} with and without applying quadratic refinement. Note that the improvement is more apparent in Test 2, which has a {\em uniformly} semi-concave solution. This agrees with a general drop in the performance of numerical schemes, especially those of higher order, in lack of uniform semi-concavity (see the discussion of this point in \cite{FF01,FF13}).

\begin{figure}[!t]
    \centering
    \includegraphics[width=0.6\textwidth]{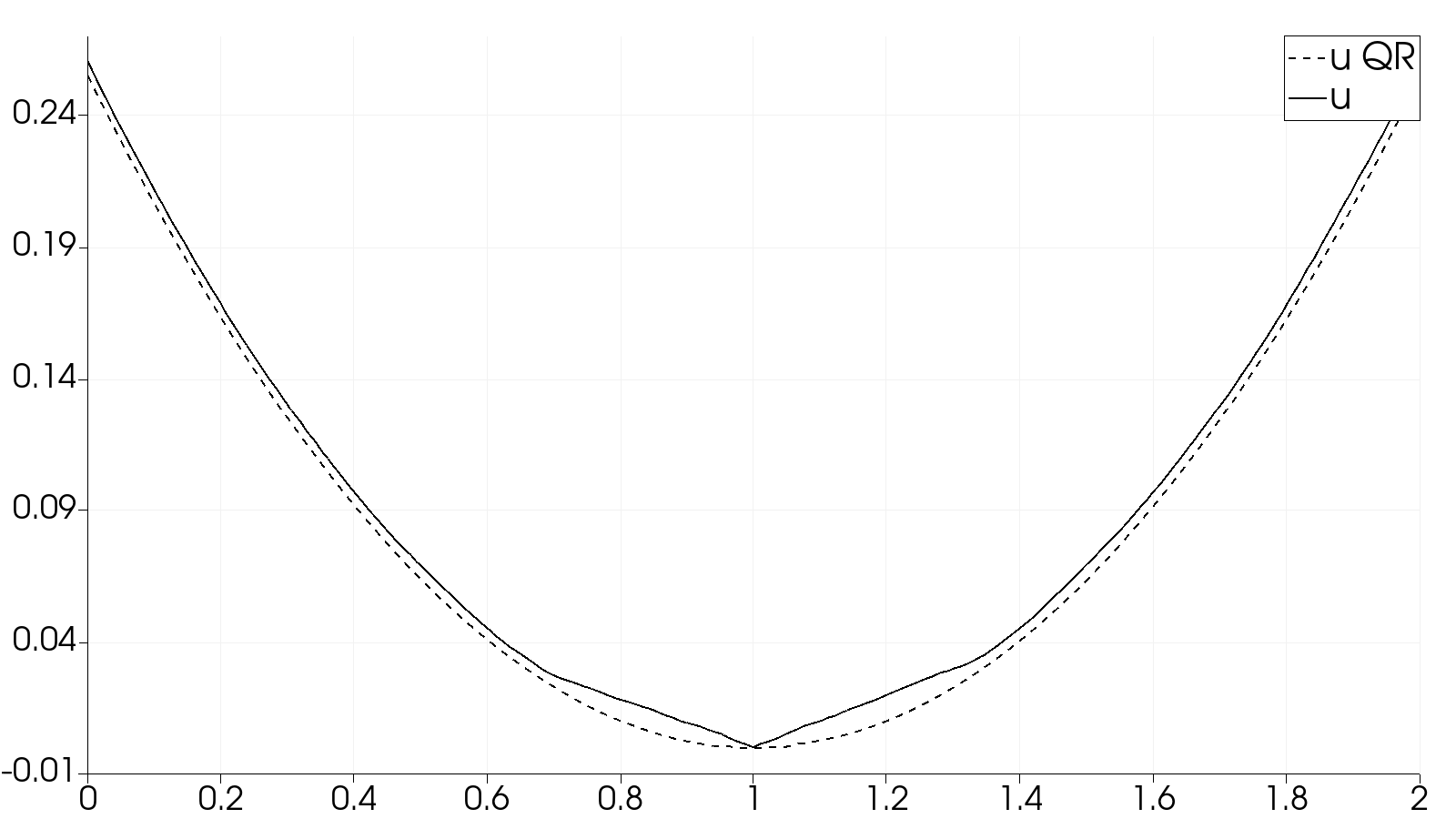}
    \caption{Sections of the numerical solutions of test $1$ at final time $T=2$. 
    The numerical solutions are indicated as $u$QR and $u$ and are obtained applying scheme \eqref{eq:hl_discr} with and without quadratic refinement, for $\Delta x= 0.0125$, $\Delta t= 0.5\sqrt{\Delta x}\approx 0.05$.}
    \label{fig:quadratic-profile}
\end{figure}

Moreover, we observe that the quadratic refinement allows us to recover the quadratic profile of the exact solution near the origin, which is otherwise lost in the previous version of the scheme.
A comparison between the solutions obtained with the two schemes is visually reported in Figure \ref{fig:quadratic-profile}, where we present the profile of both numerical solutions for test 1, where $\Omega$ is defined as a ball of radius 1, obtained sectioning along the plane parallel to the Z-axis passing through the points $(-1,0,0)$ and $(1,0,0)$.

Finally, we note that, despite significantly improving accuracy, this implementation of the quadratic refinement may cause an irregular, and lower than expected, order of convergence, even in the uniform semi-concave case. This might depend on the impossibility to enforce unisolvency (see Remark \ref{rem:acuteness}) and is worth a deeper investigation, which is, however, beyond the scope of this paper.

\subsection{Stationary case}\label{subsec:stationary}

In the second set of examples, we treat the static case, focusing in particular on the accuracy of the scheme, whereas fast solvers and parallel versions will be treated in the next subsection.

\begin{figure}[!t]
    \centering
    \includegraphics[width=.49\textwidth]{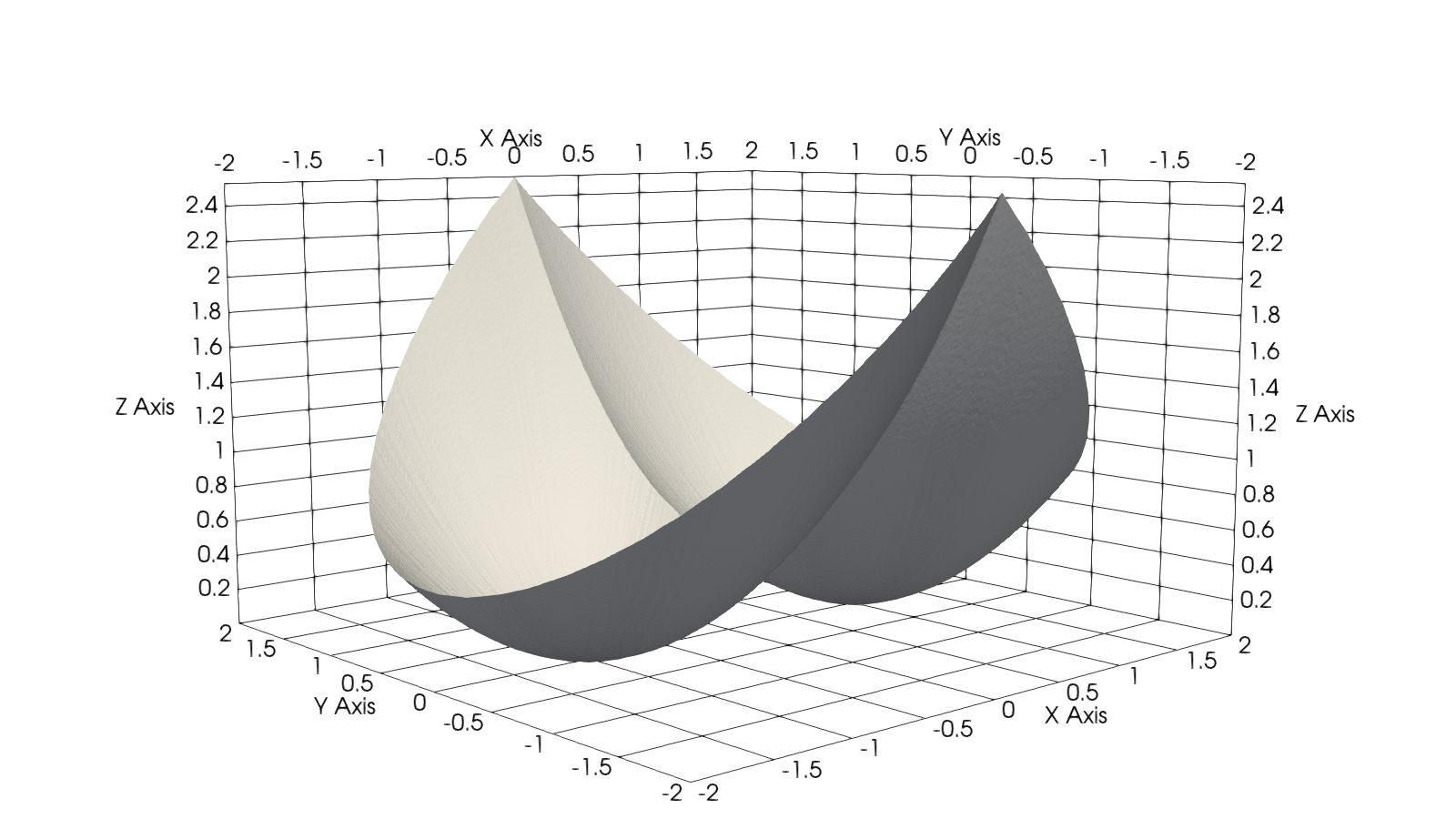}
    \includegraphics [width=.49\textwidth]{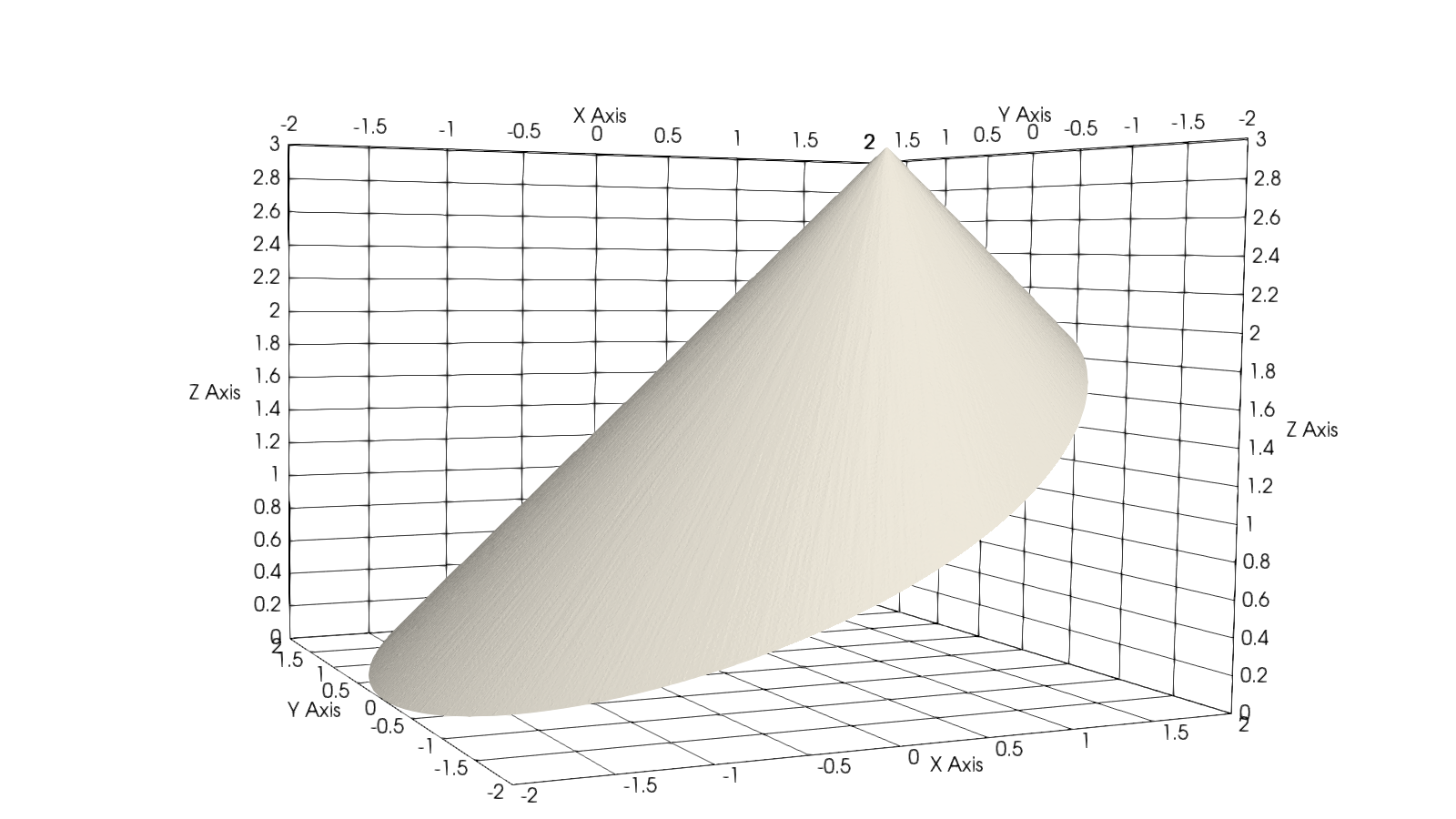}
    \caption{Tests 3--4, numerical solutions computed applying trapezoidal integration and exact policy iteration solver, for $\Delta x= 0.0125$, $\Delta t= 0.2 \sqrt{\Delta x}\approx 0.02$.}  
    \label{fig:test3-4}
\end{figure}

As for the evolutive case, we consider the quadratic Hamiltonian $H(p)= \frac 12 |p|^2$, we fix $\lambda>0$, and we study the following two tests for $x\in \Omega\subset\mathbb R^2$, where $\Omega$ is defined as the circle centered in the origin, with radius 2. In the two cases, we define the function $f$ in \eqref{eq:hj_staz} as 
\[
f_3(x)= \frac 12 (\lambda+1)\min\big\{|x-(1,0)|^2, |x+(1,0)|^2\big\}, \qquad f_4(x)= 3\lambda+\frac 1 2 - \lambda|x-(1,0)|,
\]
and enforce the boundary conditions
\[
b_3(x)\equiv 3, \qquad b_4(x)=3-\sqrt{5-2x_1}.
\]
Therefore, the associated solutions are
\[
u_3(x)= \frac 12 \min\{|x-(1,0)|^2, |x+(1,0)|^2\}, \qquad u_4(x)=3 - |x-(1,0)|,
\] 
which are shown in Fig. \ref{fig:test3-4}. Note that the boundary conditions are never active in the former example, while they are attained in the latter.

In tables \ref{tab:test3_rect}--\ref{tab:test4_rect}, we present the results obtained by applying the scheme with rectangular integration to tests 3 and 4, with the exact solutions $u_3$ and $u_4$. The value iteration starts from a constant initial guess equal to 1, and we set the fictitious time step as $\Delta t= 0.5{\Delta x}^{2/3}$, to recover a $2/3$ consistency rate. 
Tables \ref{tab:test3_trap}--\ref{tab:test4_trap} report the corresponding results for trapezoidal integration, with $\Delta t= 0.2 \sqrt{\Delta x}$, in order to obtain a unitary consistency rate.
We set $\lambda=1$, the displacement length coefficient $C=2$ and the stopping tolerance as $10^{-12}$. The experimental convergence rates are in agreement with the consistency estimates.


\begin{table}[!t]
    \centering
    \begin{tabular}{|c|c|c|c|c|c|c|c|}
        \hline
                $\Delta x$ & $\Delta t$ & $|\mathcal{X}|$ & $|\mathcal{T}|$  & $\mathcal{E}^1_{\Delta x}$ & $\text{EOC}^1$ & $\mathcal{E}^\infty_{\Delta x}$ & $\text{EOC}^\infty$ \\
                \hline\hline
                0.1 & 0.1077 & 2029 & 3928  & 0.1695 & - & 0.1318 & - \\
        \hline 
                0.05 & 0.0678 & 7919 & 15555 & 0.1125 & 0.5916 & 0.0851 & 0.6313 \\
        \hline 
                0.025 & 0.0427 & 31362 & 62142 & 0.0727 & 0.6299 & 0.0547 & 0.6381\\
        \hline  
                0.0125 & 0.0269 & 125211 & 249295  & 0.0461 & 0.6573 & 0.0348 & 0.6515\\
        \hline 
    \end{tabular}

     \caption{Test 3, experimental order of convergence in the case of rectangular integration ($\Delta t= 0.5 \Delta x^{2/3}$), and comparison of the computational performances associated to the three solvers. }
    \label{tab:test3_rect}


    \centering
    \begin{tabular}{|c|c|c|c|c|c|c|c|}
        \hline
                $\Delta x$ & $\Delta t$ & $|\mathcal{X}|$ & $|\mathcal{T}|$  & $\mathcal{E}^1_{\Delta x}$ & $\text{EOC}^1$ & $\mathcal{E}^\infty_{\Delta x}$ & $\text{EOC}^\infty$ \\
                \hline\hline
                0.1 & 0.1077 & 2318 & 4274  & 0.0588 & - & 0.0769 & - \\
        \hline 
                0.05 & 0.0678 & 8024 & 15686 & 0.0433 & 0.4777 & 0.0514 & 0.5823 \\
        \hline 
                0.025 & 0.0427 & 31698 & 62709 & 0.0268 & 0.6549 & 0.0322 & 0.6734\\
        \hline  
                0.0125 & 0.0269 & 125625 & 249877  & 0.0172 & 0.6401 &0.0205 & 0.6544\\
        \hline 
    \end{tabular}

     \caption{Test 4, experimental order of convergence in the case of rectangular integration ($\Delta t= 0.5 \Delta x^{2/3}$), and comparison of the computational performances associated to the three solvers. }
    \label{tab:test4_rect}


    \centering
    \begin{tabular}{|c|c|c|c|c|c|c|c|}
        \hline
                $\Delta x$ & $\Delta t$ & $|\mathcal{X}|$ & $|\mathcal{T}|$  & $\mathcal{E}^1_{\Delta x}$ & $\text{EOC}^1$ & $\mathcal{E}^\infty_{\Delta x}$ & $\text{EOC}^\infty$ \\
                \hline\hline
                0.1 & 0.0632 & 2029 & 3928  & 0.1599 & - & 0.0754 & - \\
        \hline 
                0.05 & 0.0447 & 7919 & 15555 & 0.0894 & 0.8384 & 0.0367 & 1.0799 \\
        \hline 
                0.025 & 0.0316 & 31362 & 62142 & 0.0471 & 0.9253 & 0.0169 & 1.0748 \\
        \hline  
                0.0125 & 0.0224 & 125211 & 249295  & 0.0233 & 1.0132 & 0.0081 & 1.0623\\
        \hline 
    \end{tabular}

     \caption{ Test 3, experimental order of convergence in the case of trapezoidal integration ($\Delta t= 0.2\sqrt{\Delta x}$) and comparison of the computational performances associated to the three solvers. }
    \label{tab:test3_trap}


    \centering
    \begin{tabular}{|c|c|c|c|c|c|c|c|}
        \hline
                $\Delta x$ & $\Delta t$ & $|\mathcal{X}|$ & $|\mathcal{T}|$  & $\mathcal{E}^1_{\Delta x}$ & $\text{EOC}^1$ & $\mathcal{E}^\infty_{\Delta x}$ & $\text{EOC}^\infty$ \\
                \hline\hline
                0.1 & 0.0632 & 2318 & 4274  & 0.0320 & - &  0.0496 & - \\
        \hline 
                0.05 & 0.0447 & 8024 & 15686 & 0.0167 & 0.9409 & 0.0190 & 1.3828 \\
        \hline 
                0.025 & 0.0316 & 31698 & 62709 & 0.0102 & 0.7076 & 0.0108 & 0.8111 \\
        \hline  
                0.0125 & 0.0224 & 125625 & 249877  & 0.0046 & 1.1400 & 0.0048 & 1.1777\\
        \hline 
    \end{tabular}

     \caption{ Test 4, experimental order of convergence in the case of trapezoidal integration ($\Delta t= 0.2\sqrt{\Delta x}$) and comparison of the computational performances associated to the three solvers. }
    \label{tab:test4_trap}

\end{table}





\subsection{Fast solvers, GPU implementation and complex geometries}
Semi-Lagrangian schemes provide a concrete example of {\em embarrassingly parallel} algorithms, meaning that the workload for the numerical solution of our Hamilton--Jacobi equations can be naturally split among the computational resources. Indeed, referring to Algorithm \ref{ALG1} for evolutive problems and to the value iteration approach for stationary problems, each node of our triangular grid can be assigned to a single thread in a parallel computation, which can perform its walks along the mesh independently, using information from the previous time step or iteration. Synchronization is needed only at the end of each iteration, but since the length of each walk is similar for each node, the amount of operations is about the same for each thread, and this clearly results in a significant speed-up of the whole computation.

Here, we implement and compare five solvers for the stationary problem: serial and parallel versions of the value iteration, serial and parallel versions of the modified policy iteration, and serial version of the exact policy iteration. The parallel codes are written in CUDA, running on a GPU server equipped with a Nvidia H100 GPU with 94Gb Ram and 14592 CUDA cores, whereas, for a fair comparison, the serial codes have been run in single-processor mode on the same platform. Note that we have given up writing an optimized parallel code for the exact policy iteration. In fact, parallel implementations of direct solvers for linear systems are typically highly problem-dependent, and with speed-ups typically not exceeding a few units (see \cite{FLRK,KOEB,N}). Moreover, from the point of view of a practical implementation, setting-up and managing performance of GPU libraries for sparse matrices in CUDA (e.g. \texttt{cusparse}) is not as simple as assembling a CUDA kernel for an iterative and explicit solver. 

In Table \ref{tab:speed-up}, we report the results obtained for test 3 of Subsection \ref{subsec:stationary}, using the same parameters and the trapezoidal quadrature in the Hopf--Lax formula, showing in particular the speed-up with respect to the serial version of the code. Note that iterative solvers achieve an impressive speed-up, even higher in the case of Value Iteration scheme. This acceleration causes iterative solvers, which are less performant in serial implementation, to outperform the Exact Policy Iteration scheme even if the latter is implemented in CUDA.


\begin{table}[!t]
    \centering
    \begin{tabular}{|c||c|c|l||c|l||c|c|l|}
        \hline
         & \multicolumn{3}{|c||}{Value Iteration} & \multicolumn{2}{|c||}{Exact PI} & \multicolumn{3}{|c|}{Modified PI} \\
        \hline\hline
            $\Delta x$ & $N$ & CPU & GPU (x) & $N$ & CPU & $N$ & CPU & GPU (x) \\
            \hline\hline
            0.1 & 394 & 1.362 & 0.066 (20x) & 12 & 0.066 & 12 & 0.149 & 0.032 (4.6x) \\
            \hline
            0.05 & 549 & 9.351 & 0.126 (74x) & 15 & 0.276 & 15 & 0.686 & 0.047 (14.6x) \\
            \hline
            0.025 & 766 & 66.347 & 0.228 (291x) & 19 & 1.676 & 19 & 4.264 & 0.077 (55x) \\
            \hline
            0.0125 & 1067 & 482.065 & 1.012 (476x) & 24 & 10.506 & 24 & 28.387 & 0.161 (176x) \\
            \hline
        \end{tabular}
    \caption{Test 3, comparison of the computational performances associated to the three solvers.}
    \label{tab:speed-up}
\end{table}

However, a key factor limiting performance is the lack of spatial locality in unstructured grids. Unlike structured meshes, where neighbouring nodes are stored contiguously in memory, an unstructured mesh requires each node to retrieve solution values and connectivity information from global memory, leading to inefficient memory access patterns. Ideally, performance could be further improved by grouping spatially close nodes into thread blocks so that neighbouring nodes share cached data, reducing redundant global memory accesses. The implementation of this strategy, however, falls once again outside the scope of the present study.  

\begin{figure}[!t]
\centering
    \includegraphics[width=.7\textwidth]{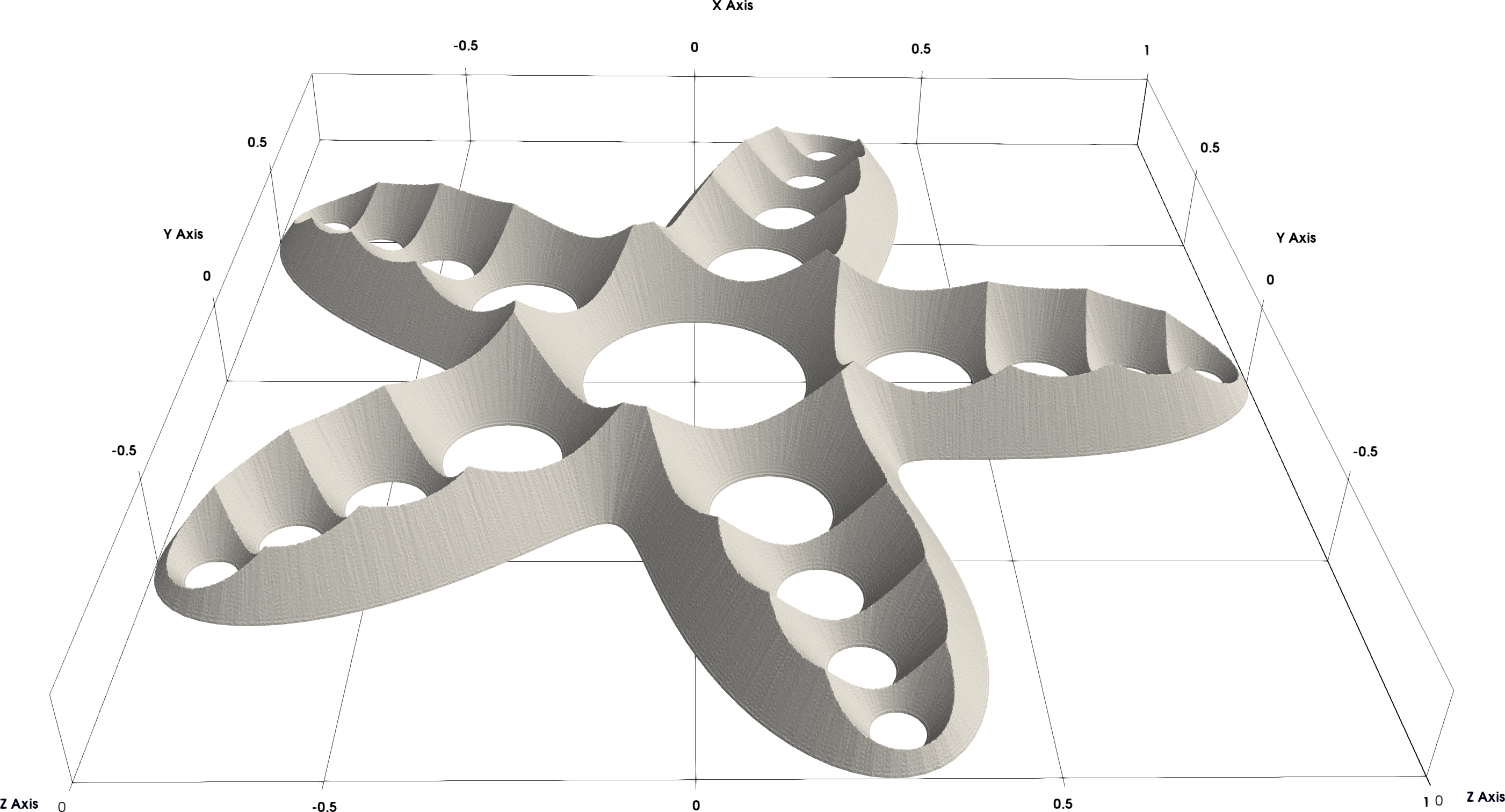}
    \caption{Numerical solution of the HJ equation \eqref{eq:hj_staz} on a flower-shaped domain with holes.}
    \label{fig:complex_geometry}
\end{figure}

We finally report the results of an experiment showing the ability of the proposed schemes to handle complex geometries in two dimensions. 
We consider a flower-shaped domain with 5 petals from which some circles of different sizes have been removed. We impose homogeneous Dirichlet conditions on all the internal and external boundaries, and we consider again the case of a stationary Hamilton--Jacobi equation with quadratic Hamiltonian. Moreover, we choose $\lambda=1$ and we set $f\equiv 1$ on the whole domain. In this case, it is easy to show that the solution of \eqref{eq:hj_staz} coincides with the Kru\v zkov transform of the distance function from the boundary of the domain \cite{FF13}. In Figure \ref{fig:complex_geometry}, we show the solution computed by the GPU version of our value iteration code. The triangular grid is composed by 330281 points and 652520 triangles, while the convergence is achieved after 37 iterations in 0.04 seconds.

\section{Conclusions}\label{sec:conclu}

In this paper, we study the unstructured implementation of a known numerical scheme for solving Hamilton--Jacobi equations using a Semi-Lagrangian approach. The key advantage of the method is that it avoids the computationally expensive interpolation step by relying solely on nodal values; our contribution here is to improve its efficiency and accuracy on unstructured space grids.

We review the theoretical accuracy and convergence analysis within the Barles–-Souganidis framework and introduce a quadratic refinement technique to improve precision while maintaining computational feasibility. For stationary problems, we compare value iteration with fast solvers of policy iteration type, discussing computational times in both serial and GPU parallel implementations, and demonstrating, in the latter case, substantial speed-up gains for large-scale problems. 

Extensive numerical tests show reliable convergence rates and robust performance even in the case of nonsmooth solutions.

\section{Aknowledgements}
S.C. is supported by the PNRR-MUR project {\em Italian Research Center on High Performance Computing, Big Data and Quantum Computing} and by the INdAM-GNCS Project code CUP E53C23001670001. R.F. is supported by the PRIN Project 2022238YY5, ``Optimal control problems: analysis, approximation and applications''.
All the Authors are members of the Gruppo Nazionale per il Calcolo Scientifico of Istituto Nazionale di Alta Matematica ``Francesco Severi'' (INdAM--GNCS).


\thebibliography{XXX}

\bibitem{BCD} M. Bardi, I. Capuzzo-Dolcetta,
Optimal control and viscosity solutions of Hamilton--Jacobi--Bellman equations. Birkhauser, Boston, 1997

\bibitem{BS} G. Barles, P.E. Souganidis, {\em Convergence of Approximation Schemes for Fully Nonlinear Second-order Equations}, Asymp. Anal., {\bf 4} (1991), 271--283.  

\bibitem{BMZ09} O. Bokanowski, S. Maroso, H. Zidani, {\em Some convergence results for Howard's algorithm}, SIAM Journal on Numerical Analysis {\bf 47} (2009), 3001--3026.

\bibitem{BR} F. Bornemann, C. Rasch, {\em Finite-element discretization of static Hamilton--Jacobi equations based on a local variational principle}, Comput. Visual. Sci. {\bf 9} (2006), 57--69.

\bibitem{BFZ} A. Bouillard, E. Faou, M. Zavidovique, {\em Fast weak–KAM integrators for separable Hamiltonian systems}, Math. Comp., {\bf 85} (2016), 85--117

\bibitem{CF} S. Cacace, R. Ferretti, {\em Efficient implementation of characteristic-based schemes on unstructured triangular grids}, Computational and Applied Mathematics, {\bf 41} (2022), 1--24

\bibitem{CT} M.G. Crandall, L. Tartar, {\em Some relations between nonexpansive and order preserving mappings}, Proc. Amer. Math. Soc., {\bf 78} (1980), 385--390

\bibitem{FF01} M. Falcone, R. Ferretti, {\em Semi-Lagrangian schemes for Hamilton--Jacobi equations, discrete representation formulae and Godunov methods}, J. of Computational Physics, {\bf 175} (2002), 559--575.

\bibitem{FF13} M. Falcone, R. Ferretti, Semi-Lagrangian approximation schemes for linear and Hamilton--Jacobi equations, SIAM, Philadelphia, 2013.

\bibitem{FF16} M. Falcone, R. Ferretti, {\em Numerical methods for Hamilton--Jacobi type equations}, Handbook of Numerical Analysis, {\bf 17}, Elsevier, 603--626.

\bibitem{FLRK} E. A. Fattah, H. Ltaief, H. Rue and D. Keyes, {\em sTiles: An Accelerated Computational Framework for Sparse Factorizations of Structured Matrices}, ISC High Performance 2025 Research Paper Proceedings (40th International Conference), Hamburg, Germany, 2025, pp. 1-14.

\bibitem{GS} J. Gianatti, F. Silva, {\em Approximation of deterministic mean field games with control-affine dynamics}, Foundations of Computational Mathematics, {\bf 24} (2002), 2017--2061.

\bibitem{KOQ} C.-Y. Kao, S. Osher, J. Qian, {\em Legendre-transform-based fast sweeping methods for static Hamilton--Jacobi equations on triangulated meshes}, J. of Computational Physics, {\bf 227} (2008), 10209--10225.

\bibitem{KOEB} M. O. Karsavuran, E. G. Ng, and B. W. Peyton. {\em GPU Accelerated Sparse Cholesky Factorization}, SC24-W: Workshops of the International Conference for High Performance Computing, Networking, Storage and Analysis. IEEE, 2024.

\bibitem{N} M. Naumov, {\em Incomplete-LU and Cholesky Preconditioned Iterative Methods Using CUSPARSE and CUBLAS}, Nvidia Whitepaper,\\\texttt{https://docs.nvidia.com/cuda/incomplete-lu-cholesky/index.html}

\bibitem{PB} M.L. Puterman, S.L. Brumelle, {\em On the convergence of policy iteration in stationary dynamic programming}, Mathematics of Operational Research {\bf 4} (1979), 60--69.

\bibitem{PS} M.L. Puterman, M.C. Shin, {\em Modified policy iteration algorithms for discounted Markov decision problems}, Management Science {\bf 24} (1978), 1127--1137.

\bibitem{SV} J. Sethian, A. Vladimirsky, {\em Ordered Upwind Methods for Static Hamilton--Jacobi Equations: Theory and Algorithms}, SIAM Journal on Numerical Analysis {\bf 41} (2003), 325--363.

\bibitem{S} C.-W. Shu, {\em High order numerical methods for time dependent Hamilton--Jacobi equations}, in ``Mathematics and computation in imaging science and information processing'' (2007), 47--91.

\bibitem{SUS}{\em SuiteSparse, A Suite of Sparse Matrix Software},\\\texttt{https://people.engr.tamu.edu/davis/suitesparse.html}.

\bibitem{TRI} {\em Triangle, A Two-Dimensional Quality Mesh Generator and Delaunay Triangulator}, \texttt{https://www.cs.cmu.edu/$\sim$quake/triangle}.
\end{document}